\documentclass[11pt]{amsart}
\usepackage{amsmath,amscd,amssymb,amsfonts,amsthm,stmaryrd,color,mathtools}
\usepackage{comment}

\usepackage{extarrows}
\usepackage[cmtip,matrix,arrow,color,poly]{xy}

\newtheorem{theorem}{Theorem}[section]
\newtheorem{corollary}[theorem]{Corollary}
\newtheorem{proposition}[theorem]{Proposition}

\newtheorem{lemma}[theorem]{Lemma}
\theoremstyle{definition}
\newtheorem{definition}[theorem]{Definition}
\theoremstyle{remark}

\newtheorem{remark}[theorem]{Remark}

\newtheorem{example}[theorem]{Example}

\newcommand{\cA}{\mathcal{A}}

\newcommand{\cD}{\mathcal{D}}

\newcommand{\cF}{\mathcal{F}}

\newcommand{\cI}{\mathcal{I}}

\newcommand{\cM}{\mathcal{M}}
\newcommand{\cN}{\mathcal{N}}

\newcommand{\cT}{\mathcal{T}}
\newcommand{\cV}{\mathcal{V}}

\newcommand{\vf}{\mathcal{T}_{\cM}}

\newcommand{\fX}{\mathfrak{X}}

\newcommand{\bN}{\mathbb{N}}
\newcommand{\bR}{\mathbb{R}}
\newcommand{\bZ}{\mathbb{Z}}

\newcommand{\severa}{\v{S}evera\ }
\newcommand{\sander}{\v{S}ander\ }
\newcommand{\Alie}{(A\to M,[\cdot,\cdot],\rho)}

\newcommand{\man}[1]   {#1\text{-}\mathcal{M}an}
\newcommand{\cob}[1]   {#1\text{-}\mathcal{C}ob}

\DeclareMathOperator{\id}{Id} \DeclareMathOperator{\im}{im}
 \DeclareMathOperator{\rk}{rank}

 \DeclareMathOperator{\End}{End}
 
\DeclareMathOperator{\VB2}{VB2}\DeclareMathOperator{\sym}{S^\bullet}

\DeclareMathOperator{\der}{Der}
\DeclareMathOperator{\Hom}{Hom}

\DeclareMathOperator{\sgn}{sgn}

\definecolor{forest}{rgb}{0.15,0.45,0.1}

\begin{document}
\sloppy
\title{A geometric characterization of $\mathbb{N}$-manifolds and the Frobenius theorem}

\author{ Henrique Bursztyn}
\address{Instituto de Matem\'atica Pura e Aplicada,
Estrada Dona Castorina 110, Rio de Janeiro, 22460-320, Brazil }
\email{henrique@impa.br}
\author{Miquel Cueca}
\address{Mathematics Institute\\Georg-August-University of G\"ottingen\\Bunsenstra{\ss}e 3-5\\G\"ottingen 37073\\Germany}
\email{miquel.cuecaten@mathematik.uni-goettingen.de}
\author{Rajan Amit Mehta}
\address{Department of Mathematics and Statistics\\  Smith College\\ 44 College Lane, Northampton, MA 01063\\ USA}
\email{rmehta@smith.edu}

\keywords{Graded manifolds, Frobenius theorem} 
\subjclass{58A50, 58A30, 53C12}

\begin{abstract}
This paper studies graded manifolds 
with local coordinates concentrated in non-negative degrees. We provide a canonical description of these objects in terms of classical geometric data and, 
building on this geometric viewpoint, 
we prove the Frobenius theorem for distributions in this graded setting. 

\end{abstract}

\maketitle
\tableofcontents

\section{Introduction}

Graded manifolds are generalizations of usual manifolds in which local coordinates carry an additional grading, most often by 
$\mathbb{Z}_2$ or $\mathbb{Z}$. They arise in several areas of mathematics and physics such as supergeometry \cite{ber:sup, fio:sup, var:book},
BV-BFV formalism and homological reduction (see \cite{cat:bvbfv} and references therein), topological field theory (e.g. AKSZ sigma models \cite{aksz, mnev:book}), double structures 
\cite{raj:tes,vor:gra}, 
derived differential geometry \cite{ber:der}, among others.  

This paper concerns $\mathbb{N}$-graded manifolds ({\em $\mathbb{N}$-manifolds} for short) of a given degree $n\in \mathbb{N}$, i.e. graded manifolds with local coordinates concentrated in degrees $0, \ldots, n$, 
and that commute in the graded sense; a precise definition is given in terms of  sheaves of graded algebras with a suitable local model, see $\S\ref{S2.1}$ (for another viewpoint, see \cite[Thm. 5.8]{rot:mon} and also \cite{gra:gro,gra:gra}). 
Such graded manifolds have become a key tool in the study of Poisson and related structures, including Lie algebroids and their higher versions 
\cite{ bon:on, sheng:ex, vai:lie}, 
as well as 
Courant algebroids and generalized geometry \cite{bur:super, roy:on, sev:some}. The goal of this paper is twofold:
we give a description of $\mathbb{N}$-manifolds by means of ordinary vector bundles and, building on this geometric viewpoint, we prove the 
Frobenius theorem for $\mathbb{N}$-manifolds, establishing the equivalence between involutive and integrable distributions.

To describe the first part of the paper,
consider an $\mathbb{N}$-manifold $\cM$ with body $M$. The local condition on the sheaf of functions $C_\cM$ implies  that the subsheaf of homogeneous
functions of degree $i$ is realized as the sheaf of sections of a vector bundle, i.e., for each $i=1,2,\ldots$ there is a vector bundle $E_{-i}\to M$ such 
that $C_\cM^i=\Gamma_{E_{-i}^*}$. Moreover, the algebra structure on $C_\cM$, $C^i_\cM\times C^j_\cM\to
C^{i+j}_\cM$,  gives rise to vector-bundle maps
\begin{equation}\label{eq:muij}
E_{-(i+j)}\to E_{-i}\otimes E_{-j}.
\end{equation}
A key observation is that, if $\cM$ is of degree $n$, then it is fully encoded in
the first $n$ vector bundles $E_{-1}, \ldots, E_{-n}$ and the maps  \eqref{eq:muij} for $i+j\leq n$. These data can be encapsulated in a coalgebra structure $\mu$ on the 
graded vector bundle ${\bf E}=\oplus_{i=1}^{n}E_{-i}$, and the coalgebra bundle $({\bf E}, \mu)$ may be seen as the ``simplest'' object from which 
the $\mathbb{N}$-manifold $\cM$ can be reconstructed.
In  Section \ref{S3:geo} we characterize the coalgebra bundles $({\bf E}, \mu)$ that arise from $\mathbb{N}$-manifolds, referred to as {\em admissible} (Def.~\ref{def-admissi}); we show, by means of a canonical construction, that there is a functor from the category of admissible $n$-coalgebra bundles 
to that of  $\mathbb{N}$-manifolds of degree $n$,
\begin{equation}\label{eq:F}
\cF: \cob{n}\to\man{n},
\end{equation}
that is an equivalence of categories, see Theorem \ref{geometrization}.
When $n=1$, this result amounts to the well-known equivalence between degree 1 $\mathbb{N}$-manifolds and ordinary vector bundles \cite{roy:on}. For $n=2$, it recovers the 
vector-bundle description of 
degree 2 $\mathbb{N}$-manifolds used in \cite[$\S$2.2]{bur:super}. These and other instances of the equivalence \eqref{eq:F} can be found in $\S$ \ref{sec:1-tan-cotan}
(see also \cite{cue:cot}).

A particular type of admissible $n$-coalgebra bundle arises from the symmetric powers of (negatively-)graded vector bundles, as described in Example \ref{splitexample}. We call those {\em split}. We show that a general admissible $n$-coalgebra bundle is always {\em noncanonically} isomorphic to one that is split; the auxiliary choices needed to establish such an isomorphism are given in Prop.~\ref{char-admis}. By verifying that the
functor \eqref{eq:F} sends split coalgebra bundles to split $\mathbb{N}$-manifolds, one recovers the fact that any $\mathbb{N}$-manifold is noncanonically isomorphic to a split one \cite[Thm.~1]{bon:on} (which can be regarded as a version of Batchelor's theorem \cite{bat:the} for $\mathbb{N}$-manifolds).

\begin{remark} 
In \cite{man:new}, Heuer and Jotz independently gave a geometric description of $\mathbb{N}$-graded manifolds of degree $n$ as certain $n$-fold vector bundles, called {\em symmetric}, extending the equivalence between degree $2$ $\mathbb{N}$-manifolds and involutive double vector bundles in \cite{lib:lac,fer:geo,mad:man} (see also \cite{grab:pol,eli:geo}). This latter result is related to the equivalence in \eqref{eq:F}, for $n=2$, as follows: admissible $2$-coalgebra bundles can be readily seen to be equivalent to special types of {\em dvb sequences} (i.e., exact sequences of vector bundles that codify double vector bundles \cite{che:ond}, see $\S$ \ref{sec:1-tan-cotan}) that correspond to involutive double vector bundles, see \cite{fer:geo}. More generally, one would expect an extension of the structural description of double vector bundles in \cite{che:ond} to $n$-fold vector bundles to clarify the direct correspondence between admissible $n$-coalgebra bundles and symmetric $n$-fold vector bundles.

\end{remark}

Towards the second goal of the paper, in Section \ref{S4} we make use of the equivalence \eqref{eq:F} to obtain a geometric description of vector fields of non-positive degrees on $\mathbb{N}$-manifolds, see Theorem \ref{thm:compatder} and $\S$ \ref{subsec:induc}. In Section
\ref{S5} we introduce distributions on $\mathbb{N}$-manifolds 
and prove in Theorem \ref{localFrob} the Frobenius theorem in this setting, asserting that a distribution is involutive if and only if it is locally spanned by coordinate vector fields. A parallel result for supermanifolds (i.e., $\mathbb{Z}_2$-graded manifolds) can be found e.g.\ in \cite{fio:sup,var:book} (see \cite{cov:fro} for an extension to
$\bZ^n_2$-graded manifolds), though the techniques involved are
rather different.

Besides it being a foundational theorem in the theory of $\bN$-graded
manifolds, our main motivation in establishing the Frobenius theorem in this context came from its role in
the study of reduction of Poisson manifolds and Courant algebroids in \cite{bur:super, cat:red}; further possible applications include, in analogy with the classical case \cite{mac:bia},  the integration of Lie-algebroid morphisms to Lie-groupoid morphisms (Lie II) in the $\mathbb{N}$-graded setting (cf. \cite{kot:glg}), that we plan to address in a separate work.

\subsection*{Acknowledgements}
H.B. and M.C. were partially supported by CNPq and Faperj.
M.C. thanks the hospitality of Smith college and IMPA that helped the development of this project. 
We thank Alejandro Cabrera, Alberto Cattaneo, Elizaveta Vishnyakova and Marco Zambon for stimulating discussions.

\section{$\mathbb{N}$-manifolds and coalgebra bundles}\label{S2:def}

In this section we introduce degree $n$ admissible coalgebra bundles, which are geometric objects defined in terms of classical vector bundles. We will show in the next section that they form a category that is equivalent to that of $\mathbb{N}$-manifolds of degree $n$.

\subsection{The category of $\mathbb{N}$-manifolds} \label{S2.1}

Given a non-negative integer $n\in\bN$, let 
$\mathbf{V}=\oplus_{i=1}^n V_{i}$ be a graded vector space with $\dim V_i=m_i$; let $\sym \mathbf{V}$ denote the graded symmetric algebra of $\mathbf{V}$.
An \emph{$\mathbb{N}$-manifold of degree $n$} (or simply an \emph{$n$-manifold})  is a ringed space $\cM=(M, C_\cM)$, where $M$ is a smooth manifold and $C_\cM$ is a sheaf of graded commutative algebras such that any point in $M$ admits a neighborhood $U$ with an isomorphism
    \begin{equation}\label{locality}
        C_\cM|_U\cong C^{\infty}_U\otimes \sym \mathbf{V},
    \end{equation}
   where the right-hand side is the sheaf of graded commutative algebras on $U$ given by the sheafification of the presheaf $U'\mapsto C^\infty(U')\otimes \sym \mathbf{V}$, for $U'\subseteq U$ open.
    We say that $\cM$ has \emph{dimension $m_0|m_1|\cdots|m_n$}, where $m_0=\dim M$.
The manifold $M$ is known as the \emph{body} of $\cM$. Global sections of $C_\cM$ are called \emph{functions} on $\cM$. We denote by $C^l_\cM \subseteq C_\cM$ the subsheaf of homogeneous sections of degree $l$. The degree of a homogeneous section $f$ is denoted by $|f|$. By the local property \eqref{locality},
\begin{equation}\label{eq:homoglocal}
C_\cM^l(U)=C^\infty(U)\otimes \mathrm{S}^l \mathbf{V}.
\end{equation}


 Let $\cM=(M,C_\cM)$ be an $n$-manifold of dimension $m_0|\cdots|m_n$. A chart $U\subseteq M$ for which  \eqref{locality} holds is called a \emph{chart} of $\cM$.  We say that 
 $$
 \{  x^\alpha, e^{\beta_i}_i\},  \quad \mathrm{ where } \; \alpha=1,\ldots, m_0,\, i=1,\ldots, n, \, \beta_i=1, \ldots, m_i, 
 $$   
 are \emph{local coordinates} of $\cM_{|U}=(U,C_\cM|_U)$  if $\{x^\alpha\}_{\alpha=1}^{m_0}$ are local coordinates on $U$ and $\{e^{\beta_i}_i\}_{\beta_i=1}^{m_i}$ form a basis of $V_i$ for all $i =1,\ldots, n$.  Then any homogeneous section of $C_\cM$ over $U$ can be expressed as a sum of functions that are smooth in $x^i$ and polynomial in $\{e^{\beta_i}_i\}$, with $|e^{\beta_i}_i|=i$.
  We  may occasionally simplify the notation and suppress the sub-indices, denoting local coordinates by 
  $\{x^\alpha, e^\beta\}$.

A morphism of $n$-manifolds $\Psi:\cM\to \cN$ is a  morphism of ringed spaces, given by a pair $\Psi=(\psi, \psi^\sharp)$, where $\psi:M \to N$ is a smooth map and $\psi^\sharp:C_\cN\to\psi_* C_\cM$ is a degree preserving morphism of sheaves of algebras over $N$. 

For each $n \in \bN$, $n$-manifolds with their morphisms form a category that we denote by $\man{n}$.

An $n$-manifold $\cM=(M, C_\cM)$ gives rise to a tower of graded manifolds
\begin{equation}\label{tower}
M=\cM_0\leftarrow\cM_1\cdots\leftarrow\cM_{n-1}\leftarrow\cM_{n}=\cM,
\end{equation}
where $\cM_r=(M, C_{\cM_r})$ is an $r$-manifold with $C_{\cM_r}$ the subsheaf of algebras of $C_\cM$ locally generated by functions of degree $\leq r$, see \cite{roy:on}.

We illustrate $n$-manifolds with some examples.

\begin{example}[Linear $n$-manifolds]\label{Ex-linear}
The most basic example of an ordinary smooth manifold of dimension $m$ is $\bR^m$.  
For the analogous picture in the graded setting, we consider a collection of nonnegative integers $m_0,\cdots,m_n$ and define the graded vector space
\begin{equation*}
 \mathbf{W}=\oplus_{i=1}^n W_i,
\end{equation*}
where $W_i=\bR^{m_i}$. The standard linear $n$-manifold of dimension $m_0|\cdots|m_n$ is defined as
\begin{equation}
    \bR^{m_0|\cdots|m_n}=\big(\bR^{m_0}, C_{\bR^{m_0|\cdots|m_n}}=C^\infty_{\bR^{m_0}}\otimes \sym \mathbf{W}\big),
\end{equation}
where $C^\infty_{\bR^{m_0}}\otimes \sym \mathbf{W}$ is the sheaf of graded algebras on $\bR^{m_0}$ obtained by sheafification of the presheaf defined on each open subset $U\subseteq \bR^{m_0}$ by $C^\infty(U)\otimes \sym \mathbf{W}$.  

 A degree $l$ homogeneous section of $C_{\bR^{m_0|\cdots|m_n}}$ over an open $U$ is an element in $C^\infty(U)\otimes \mathrm{S}^l \mathbf{W}$. An arbitrary section $f$ over $U$, once restricted to sufficiently small open subsets of $U$ (depending on $f$), are expressed by a finite sum of homogeneous sections.

 These particular $n$-manifolds are the local models for arbitrary $n$-manifolds, by \eqref{locality}.
\hfill $\diamond$
\end{example}

\begin{example}[1-manifolds]\label{1-man}
Given a vector bundle $E\to M$, we can define the $1$-manifold   $$E[1]=(M,\Gamma_{\wedge^\bullet E^*}).$$
On a chart $U$ of $M$ where $E_{|U}$ is trivializable, if we pick coordinates $\{x^\alpha\}_{\alpha=1}^m$ on $U$ and a frame $\{e^\beta\}_{\beta=1}^{\rk(E)}$  of $E^*_{|U}$ we have that $\{x^\alpha, e^\beta\}$ 
are local coordinates of $E[1]$, where $|x^\alpha|=0$ and $|e^\beta|=1$. Any $1$-manifold is canonically of this type, as we recall  in Section \ref{S3:geo}.
Two important examples of 1-manifolds are 
$$
T[1]M=(M,\Omega^\bullet_M) \qquad \text{and} \qquad T^*[1]M=(M,\fX^\bullet_M),
$$ 
where $\Omega^\bullet_M=\Gamma_{\wedge^\bullet T^*M}$ and $\fX^\bullet_M= \Gamma_{\wedge^\bullet TM}$ are the sheaves of differential forms and multivector fields on $M$, respectively. \hfill $\diamond$
\end{example}

\begin{example}[Split $n$-manifolds]\label{ex:splitnman}
The previous example can be generalized as follows. Given a non-positively graded vector bundle  $\mathbf{D}=\oplus_{i=-n+1}^{0} D_i\to M$, we define the graded manifold $\mathbf{D}[1]=(M,  \Gamma_{\sym (\mathbf{W})})$, where $\mathbf{W} = \oplus_{i=1}^{n} D_{-i+1}^*$ is the graded dual of $\mathbf{D}$ with a degree shift by $1$ (so it 
is a graded vector bundle concentrated in degrees $1$ to $n$), $\sym$ denotes, as before, the graded symmetric product, and $\Gamma_{\sym (\mathbf{W})}$ is the sheaf of graded algebras given by  $\oplus_{l=0}^\infty \Gamma_{\mathrm{S}^l (\mathbf{W})}$. Such $n$-manifolds are known as {\em split}, see e.g. \cite{bon:on}. \hfill $\diamond$
\end{example}

\begin{remark}[Graded duals]\label{rem:duals} Given a graded vector space $\mathbf{V}=\oplus_{i\in\bZ} V_i$, its dual is the graded vector space $\mathbf{V}^*$ defined by $(\mathbf{V}^*)_j = (V_{-j})^*$.
\hfill $\diamond$
\end{remark}

\begin{example}[Shifted tangent and cotangent bundles of $n$-manifolds] \label{shift-tan-cotan}
One can generalize the shifted tangent and cotangent bundles of ordinary manifolds considered in
Example \ref{1-man} to arbitrary $n$-manifolds $\cM=(M, C_\cM)$, see e.g. \cite{cue:cot, raj:tes}. We will give a description using charts and coordinates (cf. \cite{roy:on}), so that $T[k]\cM$ and $T^*[k]\cM$ are defined by gluing local sheaves.

Let $U$ be a chart of $\cM$  with coordinates $\{x^\alpha, e^{\beta}\}$,
so that $\{x^\alpha\}$ are coordinates on $U$  
and $C_\cM(U)=C^\infty_M(U)\otimes \sym\mathbf{V}$, with $\{e^\beta\}$ a basis of $\mathbf{V}$. 
Just as for ordinary smooth manifolds, we can introduce the differentials of the coordinates, $\{dx^\alpha, de^\beta\}$, with the transformation rule 
\begin{equation*}
    \left\{\begin{array}{ll}
        \widehat{x}^\alpha=  & F^\alpha(x), \\
        \widehat{e}^\beta=  & F^\beta(x, e)
    \end{array} \right.\quad \Longrightarrow \quad\left\{\begin{array}{ll}
        d\widehat{x}^\alpha=  & \frac{\partial F^\alpha(x)}{\partial x^{\alpha'}}dx^{\alpha'}, \\
        d\widehat{e}^\beta= & \frac{\partial F^\beta(x, e)}{\partial x^\alpha}dx^\alpha+\frac{\partial F^\beta(x, e)}{\partial e^{\beta'}}de^{\beta'}. 
    \end{array}\right.
\end{equation*}
Dualizing we can also define $ \{\frac{\partial}{\partial x^\alpha}, \frac{\partial}{\partial e^\beta}\}$, with 
the corresponding transformation rule.

For $k\geq 1$ we define the $(n+k)$-manifold $T[k]\cM=(M, C_{T[k]\cM})$ by
$$ 
C_{T[k]\cM}(U)= C_M^\infty(U)\otimes \sym\mathbf{W}\quad\text{with}\quad \mathbf{W}=\mathbf{V}\oplus\langle dx^\alpha, de^\beta \rangle,
$$
with $|dx^\alpha|=k$,  $|de^\beta|=|e^\beta|+k$.
Similarly, for $k\geq n+1$ (to ensure that all  coordinates are in positive degrees), we define the $k$-manifold $T^*[k]\cM=(M, C_{T^*[k]\cM})$ by
$$ 
C_{T^*[k]\cM}(U)= C_M^\infty(U)\otimes \sym\mathbf{W}\quad\text{with}\quad \mathbf{W}=\mathbf{V}\oplus  \left \langle  \frac{\partial}{\partial x^\alpha}, \frac{\partial}{\partial e^\beta} \right \rangle,
$$ 
with $\Big |\frac{\partial}{\partial x^\alpha}\Big |=k$, $\Big |\frac{\partial}{\partial e^\beta}\Big |=-|e^\beta|+k$. 
\hfill $\diamond$
\end{example}

\begin{example}[Cartesian product of $n$-manifolds]\label{car-pro}
Given $n$-manifolds $\cM_i=(M_i,C_{\cM_i})$, $i=1,2$, we define an $n$-manifold 
$\cM_1\times\cM_2=(M_1\times M_2, C_{\cM_1\times\cM_2})$, where
\begin{equation*}
C_{\cM_1\times\cM_2}=C_{\cM_1}\widehat{\otimes}\ C_{\cM_2},
\end{equation*}
and the hat denotes the usual completion on the product topology,  see e.g. \cite{fio:sup}. Concretely, on open rectangles $U_1\times U_2$, with $U_i\subset M_i$ satisfying $C_{\cM_i}(U_i)\cong C^\infty_{M_i}(U_i)\otimes \sym \mathbf{V}_i$ we have
\begin{equation}\label{prod}
C_{\cM_1\times\cM_2}(U_1\times U_2)\cong C^\infty_{M_1\times M_2}(U_1\times U_2)\otimes \sym (\mathbf{V}_1\oplus \mathbf{V}_2),
\end{equation}
where $\mathbf{V}_1\oplus \mathbf{V}_2$ denotes the usual graded direct sum. Since open rectangles form a basis of the product topology, \eqref{prod} is sufficient to define $C_{\cM_1 \times \cM_2}$. 
\hfill $\diamond$
\end{example}

\subsection{Admissible coalgebra bundles} \label{S2.2}
Any $n$-manifold $\cM= (M, C_\cM)$ codifies classical geometric information as follows.
The multiplication on $C_\cM$  satisfies
\begin{equation}\label{mult}
C^i_\cM\cdot C^j_\cM\subseteq C^{i+j}_\cM,
\end{equation}
and since $C^0_\cM= C^\infty_M$ (as a consequence of \eqref{locality}), each $C^i_\cM$ is a $C^\infty_M$-module, which is locally free and finitely generated (see \eqref{eq:homoglocal}). It follows that there exist vector bundles $E_{-i}\to M$ such that 
$$
C^i_\cM\cong \Gamma_{E_{-i}^*}.
$$
It turns out that the fact that $\cM$ is of degree $n$ implies that $C_\cM$ is completely determined by information involving just the vector bundles $E_{-1},\cdots, E_{-n}$. More precisely,
since the multiplication in \eqref{mult} is $C^\infty_M$-bilinear, it 
defines vector bundle maps $E_{-(i+j)}\to E_{-i}\otimes E_{-j}$, and we will see (in Lemma~\ref{lem:reconst}) that $\cM$ can be reconstructed from such maps for $i+j\leq n$.

Given vector bundles $E_{-1}, \ldots, E_{-n}$ and vector bundle maps $ E_{-(i+j)}\to E_{-i}\otimes E_{-j}$ (covering the identity map on $M$), one could ask whether this data arises from an $n$-manifold. This will lead us to the notion of {\em admissible coalgebra bundles}.

\begin{definition}
    Given $n\in\bN$, an \emph{$n$-coalgebra bundle} is a pair $(\mathbf{E}\to M, \mu)$ where $\mathbf{E}=\oplus_{i=-n}^{-1} E_i\to M$ is a negatively graded finite-dimensional vector bundle and $\mu: \mathbf{E}\to \mathbf{E}\otimes \mathbf{E}$ is a degree-preserving vector bundle map (covering the identity map on $M$) that is coassociative and (graded) cocommutative, i.e. 
$$
(\mathrm{Id}\otimes \mu)\circ \mu=(\mu\otimes \mathrm{Id})\circ \mu,\qquad \mu=\tau\cdot \mu,
$$
where $\tau \in S_2$ acts via the braiding map that exchanges the two components with the appropriate sign, i.e., $\tau(e\otimes e')=(-1)^{i+j} e'\otimes e$, for $e\in E_i, \ e'\in E_j$. 
\end{definition}
A morphism  of $n$-coalgebra bundles, $(\mathbf{E}\to M,\mu) \to (\mathbf{F}\to N,\nu)$, 
is a degree-preserving vector-bundle map $\Phi: \mathbf{E}\to \mathbf{F}$ satisfying $\nu\circ \Phi = (\Phi\otimes \Phi)\circ \mu$.

The coassociative and cocommutative assumptions put constraints on the image of $\mu$, that we describe more explicitly. For $k \geq 0$, let 
\[\mu^k: \mathbf{E} \to  \otimes^{k+1}  \mathbf{E}\]
be defined by $\mu^k = (\id \otimes \cdots \otimes \id \otimes \mu) \circ \cdots \circ (\id \otimes \mu) \circ \mu\circ \id$. 
Then coassociativity and cocommutativity imply that, for all $e \in \im(\mu)$, the identity
\begin{equation} \label{coass}
\tau \cdot (\mu^k \otimes \mu^\ell)(e) = \tau' \cdot (\mu^{k'} \otimes \mu^{\ell'})(e)
\end{equation}
holds for all $k$, $k'$, $\ell$, $\ell'\in\mathbb{Z}_{\geq 0}$ such that $k + \ell = k' + \ell'$, and for all permutations $\tau, \tau' \in S_{k+\ell+2}$, acting via the (graded) braiding action\footnote{For a graded vector space
$\mathbf{V}$, the (graded) braiding action of 
the symmetric group $S_k$ on $\otimes^k \mathbf{V}$  is $\tau\cdot (v_1\otimes\cdots\otimes v_k)=\epsilon(\tau)\ v_{\tau(1)}\otimes\cdots\otimes v_{\tau(k)}$,  where $\epsilon(\tau)$ is the Koszul sign arising from the supermathematics convention that a minus sign appears whenever two consecutive odd elements are swapped.}.

\begin{definition}
We denote by $\mathbf{K}^\mu \subseteq \mathbf{E} \otimes \mathbf{E}$ the space of all $e \in \mathbf{E} \otimes \mathbf{E}$ satisfying \eqref{coass}. The subspace of degree $i$ elements of $\mathbf{K}^\mu$ is denoted by $K^\mu_i$. 
\end{definition}

As noted above, coassociativity and cocommutativity of $\mu$ imply that $\im(\mu) \subseteq \mathbf{K}^\mu$.

\begin{definition}\label{def-admissi}
An $n$-coalgebra bundle $(\mathbf{E}\to M, \mu)$ is \emph{admissible} if, $\forall \, i= -n, \ldots, -2$,  $$ 
\im(\mu_i) = {K}^\mu_i, 
$$
where 
$$
\mu_i=\mu|_{E_i}.
$$ 
(Note that the condition is trivial for $i=-1$.)
For each $n\in\bN$, we denote by $\cob{n}$ the category of admissible $n$-coalgebra bundles.
\end{definition}

\begin{remark}\label{admissi-rmk}
Since $\mathbf{E}\to M$ is concentrated in negative degrees, we see that the definition of ${K}^\mu_i$ via  \eqref{coass} only involves $\mu_l$ for $l>i$.
So the admissibility condition at level $i$ relates $\mu_i$ to the components $\mu_l$ for $l>i$. \hfill $\diamond$
\end{remark}

We illustrate the previous definition with some examples.

\begin{example}[1-coalgebra bundles] 
\label{cob1}
Let $\mathbf{E}\to M$ be a graded vector bundle concentrated in degree $-1$, $\mathbf{E}=E_{-1}$. 
Then $\mathbf{E} \otimes \mathbf{E} = E_{-1} \otimes E_{-1}$ is concentrated in degree $-2$, so any coalgebra-bundle structure must be trivial, and hence admissible. It follows that $\cob{1}$ coincides with the category $\cV ect$ of vector bundles.
\hfill $\diamond$
\end{example}

\begin{example}[2-coalgebra bundles]
\label{cob2}
For a 2-coalgebra bundle $\mathbf{E}= E_{-2}\oplus E_{-1}$, the coalgebra structure is determined by $\mu_{-2}: E_{-2}\to E_{-1}\otimes E_{-1}$. For the admissibility condition we need to determine $K^\mu_{-2}$. Note that, in the case at hand, equation \eqref{coass} is determined by $k=l=k'=l'=0, \ \tau=(12)$ and $\tau'=id$, showing that $K_{-2}^\mu = E_{-1}\wedge E_{-1}$ (viewed as a subspace of $E_{-1}\otimes E_{-1}$ in the usual way). So the admissibility condition amounts to
$$
\im(\mu_{-2}) = E_{-1}\wedge E_{-1}.
$$
So any $2$-coalgebra bundle is determined by two vector bundles $E_{-1}, E_{-2}$ and a map $\mu:E_{-2}\to E_{-1}\otimes E_{-1}$ satisfying the previous condition.
The category $\cob{2}$ is thus equivalent to the category whose objects are triples 
$$(E_{-1}\to M, E_{-2}\to M, \mu:E_{-2}\twoheadrightarrow E_{-1}\wedge E_{-1}),
$$ 
where $\mu$ is a surjective vector-bundle map over the identity on $M$, and 
 morphisms 
 $$
 (E_{-1}, E_{-2}, \mu) \longrightarrow  (F_{-1}, F_{-2}, \nu)
 $$
 are triples $(\phi,\Phi_{-1},\Phi_{-2})$, where $(\phi,\Phi_i):(E_i\to M)\to (F_i\to N)$, $i=-1,-2$,  are vector-bundle maps satisfying $\nu\circ \Phi_{-2}=(\Phi_{-1}\wedge\Phi_{-1})\circ\mu$. This category was introduced in \cite[$\S$ 2.2]{bur:super}, where it is called $\VB2$. \hfill $\diamond$
\end{example}

\begin{example}[3-coalgebra bundles]
\label{cob3}
On a 3-coalgebra bundle $(\mathbf{E},\mu)$, notice that
\begin{align*}
K^{\mu}_{-3}&=\{ e\otimes p\in E_{-1}\otimes E_{-2}\ |\ e\otimes\mu_{-2}(p)=\sgn(\tau)\tau (e\otimes \mu_{-2}(p))\quad \text{for }\tau\in S_3\}\\ &= (\mathrm{Id}_{E_{-1}}\otimes \mu_{-2})^{-1}(\wedge^3 E_{-1}).
\end{align*}
Building on the previous example,
 an admissible $3$-coalgebra bundle is determined by three vector bundles $E_{i}\to M$, $i = -1, -2, -3$, and two vector-bundle morphisms 
 $$ 
 \mu_{-2}:E_{-2}\to E_{-1}\otimes E_{-1}\quad\text{and}\quad \mu_{-3}:E_{-3}\to E_{-1}\otimes E_{-2}
 $$
satisfying 
\begin{equation*}
    (\mu_{-2}\otimes\id)\circ\mu_{-3}=(\id\otimes\mu_{-2})\circ\mu_{-3},\quad \im(\mu_{-2})=E_{-1}\wedge E_{-1}, \quad  \im(\mu_{-3})= (\mathrm{Id}_{E_{-1}}\otimes \mu_{-2})^{-1}(\wedge^3 E_{-1}). 
\end{equation*}
A morphism of $3$-coalgebra bundles from $(\mathbf{E}\to M,\mu)$ to $(\mathbf{F}\to N, \nu)$ amounts to three vector-bundle morphisms $\Phi_i:E_i \to F_i$ covering the same base map $\phi:M\to N$ satisfying
$$
\nu_{-2}\circ\Phi_{-2}=(\Phi_{-1}\otimes\Phi_{-1})\circ\mu_{-2}\quad\text{and}\quad \nu_{-3}\circ\Phi_{-3}=(\Phi_{-1}\otimes\Phi_{-2})\circ\mu_{-3}.
$$
The category $\cob{3}$ has objects $(E_{-1}, E_{-2}, E_{-3}, \mu_{-1}, \mu_{-2})$ and morphisms $(\Phi_{-1}, \Phi_{-2}, \Phi_{-3})$ as just described. 
\hfill $\diamond$
\end{example}

\begin{example}[Truncation]\label{extrunc}
Let $(\mathbf{E}\to M,\mu)$ be an admissible $n$-coalgebra bundle. For each $k=0,\ldots,n$, the $k$-truncation $(\mathbf{E}^{\leq k}\to M, \mu^{\leq k})$ is the $k$-coalgebra bundle given by  
$$(\mathbf{E}^{\leq k})_i=\left\{\begin{array}{ll}
    E_i & \text{ if }  -k\leq i \leq -1 \\
    0 & \text{ otherwise,}
\end{array}\right.\quad \mathrm{ and }\qquad \mu^{\leq k}=\mu_{|\mathbf{E}^{\leq k}}.
$$
Hence $\mu^{\leq k}_i = \mu_i$, for  $-k\leq i \leq -1$. Note also that  
\begin{equation}\label{eq:Ktrunc}
K_i^\mu = K_i^{\mu^{\leq k}}, \qquad \mbox{ for } i\geq -k-1,
\end{equation}
as a consequence of  Remark \ref{admissi-rmk}. Since $\mathbf{E}$ is admissible, for $i\geq -k$ we have that
$$
\im(\mu_i^{\leq k}) = \im(\mu_i) =  K_i^\mu =  K_i^{\mu^{\leq k}}, 
$$
so the truncation $\mathbf{E}^{\leq k}$ is also admissible.
\hfill $\diamond$
\end{example}

\begin{example}[Split coalgebra bundles]\label{splitexample}
 Given a graded vector bundle    $\mathbf{D} = \oplus_{i=-n}^{-1} D_i$ we define an admissible $n$-coalgebra bundle $(\mathbf{E}_\mathbf{D}, \mu^s)$ by
\begin{equation*}
    (\mathbf{E}_\mathbf{D})_i=
 \left(\sym \mathbf{D}\right)_i \quad \text{for } -n\leq i\leq -1,
\end{equation*}
 where $(\sym \mathbf{D})_i$ denotes the subspace of degree $i$ elements in the graded symmetric product of $\mathbf{D}$, and $\mu^s:\mathbf{E}_\mathbf{D}\to\mathbf{E}_\mathbf{D}\otimes \mathbf{E}_\mathbf{D}$ is given in components by
\begin{equation*}
    \langle\mu^s_i(e),\eta\otimes\nu\rangle=\langle e, \eta\cdot\nu\rangle
\end{equation*}
where $e\in (\mathbf{E}_\mathbf{D})_i,\ \eta\in(\mathbf{E}_\mathbf{D})^*_j=(\sym \mathbf{D}^*)_{-j},\ \nu\in(\mathbf{E}_\mathbf{D})^*_{i-j}=(\sym \mathbf{D}^*)_{j-i}$ with $-n+1\leq j\leq-1$ (see Remark~\ref{rem:duals}) and $``\cdot"$ is the graded symmetric product
$$
(\sym \mathbf{D}^*)_{-j}\otimes (\sym \mathbf{D}^*)_{j-i}\to (\sym \mathbf{D}^*)_{-i}.
$$ 
The fact that $(\mathbf{E}_\mathbf{D},\mu^s)$ is admissible follows directly from the identities 
\begin{equation}\label{eq:Kmu}
(\mathbf{E}_\mathbf{D})_i=D_i\oplus (\sym \mathbf{D}^{\leq-i-1})_i\quad \text{and}\quad K^{\mu^s}_i=(\sym \mathbf{D}^{\leq-i-1})_i\quad \text{for}\ -n\leq i\leq -2,
\end{equation}
where $\mathbf{D}^{\leq-i-1}=\oplus_{j=i+1}^{-1}D_j.$  In particular we also have that $D_i=\ker \mu^s_i$ for all $-n\leq i\leq -1.$
Such $n$-coalgebra bundles of type $(\mathbf{E}_\mathbf{D}, \mu^s)$ are called {\em split}.
\hfill $\diamond$
\end{example}

It follows from the previous example that if a coalgebra bundle is isomorphic to one that is split, then it is admissible. Conversely, we now show that any admissible $n$-coalgebra bundle is isomorphic to a split one.

\begin{proposition}[Admissible coalgebra bundles are split]\label{char-admis}
Let $(\mathbf{E}\to M,\mu)$ be an admissible $n$-coalgebra bundle. Then, for each $i=-n,\ldots,-1$, we have a short exact sequence of vector bundles
$$
0\longrightarrow \ker(\mu_i) \longrightarrow E_i \longrightarrow K^\mu_i \longrightarrow 0,
$$ 
and a choice of splittings of these sequences gives rise to an identification
$$
(\mathbf{E},\mu) {\cong} (\mathbf{E}_\mathbf{D},\mu^s),
$$
where  $(\mathbf{E}_\mathbf{D},\mu^s)$ is the split $n$-coalgebra bundle associated with the graded vector bundle $\mathbf{D}:=\oplus_{i=-n}^{-1} \ker\mu_i$ (as in Example~\ref{splitexample}).
\end{proposition}

\begin{proof}
The proof will be by induction on $n$. For $n=1$ the proposition holds since the coalgebra structures are trivial and $E_{-1}=D_{-1}$ (see Example~\ref{cob1}). 
Now assume that the proposition holds for admissible $(n-1)$-coalgebra bundles. 

Given  an admissible $n$-coalgebra bundle $(\mathbf{E},\mu)$,
 $(\mathbf{E}^{\leq n-1}, \mu^{\leq n-1})$ is an admissible $(n-1)$-coalgebra bundle (see Example \ref{extrunc}). By the induction hypothesis we have an isomorphism 
$$ 
\quad\Phi^{\leq n-1}:(\mathbf{E}^{\leq n-1}, \mu^{\leq n-1})\to (\mathbf{E}_{\mathbf{D}^{\leq n-1}},\mu^s)$$
where $\mathbf{D}^{\leq n-1}= \oplus_{i=-n+1}^{-1} \ker(\mu_i) $. 
We can now use  the isomorphism  $\Phi^{\leq n-1}$ to conclude that 
$$
\im(\mu_{-n})=(K^{\mu}_{-n}) = (K^{\mu^{\leq n-1}}_{-n})\cong (K^{\mu^s}_{-n})=(\sym \mathbf{D}^{\leq n-1})_{-n},
$$
where the first equality 
is the admissibility condition at level $-n$ for $\mathbf{E}$ and the second follows from \eqref{eq:Ktrunc}. Therefore $E_{-n}$ and $(\mathbf{E}_\mathbf{D})_{-n}$ fit into the following diagram 
\begin{equation}\label{ex-En}
\begin{array}{c}
 \xymatrix{0\ar[r]& \ker(\mu_{-n})\ar[r]\ar[d]^{\id}&  E_{-n}\ar[r]^{\mu_{-n}}& K^{\mu}_{-n}\ar[r]\ar[d]_{\cong}^{\Phi^{\leq n-1}}& 0\\
0\ar[r]& \ker(\mu_{-n}) \ar[r]&(\mathbf{E}_\mathbf{D})_{-n}\ar[r]^{\mu^s_{-n}}& (\sym \mathbf{D}^{\leq n-1})_{-n}\ar[r]& 0.}
\end{array}
\end{equation}
A splitting of the top exact sequence of vector bundles identifies $E_{-n}$ with $(\mathbf{E}_\mathbf{D})_{-n}$ and extends $\Phi^{\leq n-1}$ to a coalgebra bundle isomorphism 
 $(\mathbf{E},\mu) {\stackrel{\sim}{\rightarrow}} (\mathbf{E}_\mathbf{D},\mu^s)$.
\end{proof}


\section{The equivalence}\label{S3:geo}

For each $n\in\bN$, we define in this section a functor that establishes an equivalence of categories from $\cob{n}$ to $\man{n}$. This {\em geometrization functor} thus provides a description  of $n$-manifolds by (admissible) $n$-coalgebra bundles, which are ordinary differential-geometric objects.
It extends the known correspondence between $1$-manifolds and vector bundles, as well as the characterization of $2$-manifolds by the $\VB2$-category (see Example~\ref{cob2}) in \cite[$\S$ 2.2]{bur:super}. This functor has also been used 
 in \cite{cue:cot} to study graded cotangent bundles $T^*[n]A[1]$.

\subsection{The functor}\label{S3.1}

Let $(\mathbf{E}\to M,\mu)$ be an $n$-coalgebra bundle, and let $\mu^*:\Gamma_{\mathbf{E}^*}\times \Gamma_{\mathbf{E}^*}\to \Gamma_{\mathbf{E}^*}$ denote the dual algebra structure. Recall that
$\mathbf{E}^*=\oplus_{i=1}^n E_{-i}^*$.

In the sheaf of algebras  $( \Gamma_{\sym \mathbf{E}^*},\cdot)$ we define the following sheaf of ideals: for an open subset $U\subseteq M$,
\begin{equation*}
\cI_\mu(U)=\langle \omega\cdot\eta - \mu^*(\omega,\eta) \rangle, 
\end{equation*}
where $\omega,\eta\in\Gamma_{\mathbf{E}^*}(U)$, $|\omega|+|\eta|\leq n$, and
$``\cdot"$ denotes the graded symmetric product.

\begin{remark}\label{lowdeg}
Note that $\cI_\mu$ is a sheaf of homogeneous ideals with respect to the grading of $ \mathbf{E}^*$. Therefore $ \Gamma_{\sym \mathbf{E}^*}/\cI_\mu$ is also a  sheaf of graded algebras. Note that, for $1\leq i\leq n$,
\begin{equation*}
( \Gamma_{\sym \mathbf{E}^*}/\cI_\mu)_i=\Gamma_{E_{-i}^*}.
\end{equation*}
\hfill $\diamond$
\end{remark}

An  $n$-coalgebra bundle  $(\mathbf{E},\mu)$ gives rise to a ringed space by 
\begin{equation}\label{eq:ringed}
    (M,  \Gamma_{\sym\mathbf{E}^*}/\cI_\mu).
\end{equation}
Given a coalgebra morphism $\Phi:(\mathbf{E}\to M,\mu)\to (\mathbf{F}\to N,\nu)$, with $\phi=\Phi_{|_M}:M\to N$, we also denote  by  $\Phi:\sym \mathbf{E}\to \sym \mathbf{F}$ its natural extension to the symmetric products.  The fact that $\Phi$ is a coalgebra map implies that $\Phi^*:\Gamma_{\mathbf{F}^*}\to \phi_*\Gamma_{\mathbf{E}^*}$ is a map of sheaves of algebras, so $\Phi^* \cI_\nu\subseteq \cI_\mu$. Thus, $\Phi^*$ induces a well-defined map on the quotients by the ideals, denoted by $\phi^\sharp$, in such a way that
\begin{equation}\label{eqMor}
    (\phi,\phi^\sharp):(M, \Gamma_{\sym\mathbf{E}^*}/\cI_\mu)\to (N, \Gamma_{\sym\mathbf{F}^*}/\cI_\nu)
\end{equation}
is a morphism of  ringed spaces.

We need the following result to show that, when $(\mathbf{E},\mu)$ is admissible, the ringed space \eqref{eq:ringed} is an $n$-manifold.

 \begin{lemma}\label{lemma2}
Let  $\mathbf{D}=\oplus_{i=-n}^{-1}D_{i}\to M$ be a graded vector bundle. Consider the  coalgebra bundle $(\mathbf{E}_\mathbf{D}, \mu^s)$ of Example \ref{splitexample}. Then
\begin{equation*}
     \Gamma_{\sym\mathbf{E}_\mathbf{D}^*}/\mathcal{I}_{\mu^s}=\Gamma_{\sym  \mathbf{D}^*}.
\end{equation*}
\end{lemma}

\begin{proof}
Denote the symmetric product on $\Gamma_{\sym \mathbf{E^*_D}}$ by $``\cdot"$ and the symmetric product on $\Gamma_{\sym\mathbf{D}^*}$ by $``\vee"$. 
Since $(\mathbf{E}^*_\mathbf{D})_i=(\sym \mathbf{D}^*)_i$ for $i=1,\ldots,n$, we can view $\mathbf{D}^*\subseteq \mathbf{E}^*_\mathbf{D}$, so we have a natural inclusion 
$\Gamma_{\sym  \mathbf{D}^*} \hookrightarrow \Gamma_{\sym\mathbf{E}_\mathbf{D}^*}$, given by $d_1\vee \ldots \vee d_k \mapsto d_1\cdot\ldots\cdot d_k$. Consider the quotient homomorphism $\Gamma_{\sym\mathbf{E}_\mathbf{D}^*}\twoheadrightarrow (\Gamma_{\sym\mathbf{E}_\mathbf{D}^*}/\mathcal{I}_{\mu^s})$, $\xi \mapsto \overline{\xi}$.
From the definition of $\mathcal{I}_{\mu^s}$, one can readily see that the composite map 
\begin{equation}\label{eq:symmap}
\Gamma_{\sym  \mathbf{D}^*}\rightarrow \Gamma_{\sym\mathbf{E}_\mathbf{D}^*}/\mathcal{I}_{\mu^s}, \quad {d}_1\vee \ldots \vee {d}_k \mapsto \overline{d_1}\cdot\ldots\cdot \overline{d_k}
\end{equation}
is still an injection. To see that it is also a surjection, recall that 
$\mu^{s*}(\xi,\eta)=\xi\vee \eta$ (see Example~\ref{splitexample}), so the image of  a local section
$\xi = {d}_1\vee \ldots \vee {d}_l$ of $\Gamma_{(\mathbf{E}^*_\mathbf{D})_i}=\Gamma_{(\sym \mathbf{D}^*)_i}$
under the quotient by $\mathcal{I}_{\mu^s}$ is $\overline{\xi} = \overline{d_1}\cdot\ldots\cdot \overline{d_l}$. It follows that an arbitrary local section in $\Gamma_{\sym\mathbf{E}_\mathbf{D}^*}/\mathcal{I}_{\mu^s}$ is a symmetric expression of elements $\overline{d_j}$ in $\Gamma_{\mathbf{D}^*}$, showing that \eqref{eq:symmap} is surjective.
\end{proof}

\begin{proposition}\label{F-well-def}
For each non-negative integer $n$ the functor $\cF: \cob{n}\to\man{n}$ given  by
\begin{equation*}
\cF(\mathbf{E}\to M, \mu)=(M,  \Gamma_{\sym\mathbf{E}^*}/\cI_\mu)\qquad\text{and}\qquad
\cF(\Phi)=(\phi,\phi^\sharp):\cF(\mathbf{E}, \mu)\to \cF(\mathbf{F}, \nu)
\end{equation*}
as in \eqref{eqMor} is well defined. 
\end{proposition}
\begin{proof}
We start by showing that $\cF$ is well defined at the level of objects. Assuming that the coalgebra bundle $(\mathbf{E}\to M, \mu)$ is admissible we must prove that 
\begin{equation}\label{eq:M}
\cM=\cF(\mathbf{E}\to M, \mu)=(M,C_\cM= \Gamma_{\sym \mathbf{E}^*}/\cI_\mu)
\end{equation}
satisfies condition \eqref{locality}. By the admissibility condition,  Prop.~ \ref{char-admis} and Lemma~\ref{lemma2} imply that $\Gamma_{\sym \mathbf{E}^*}/\cI_\mu \cong \Gamma_{\sym \mathbf{D}^*}$, for $\mathbf{D}^*=\oplus_{i=1}^n (\ker \mu_{-i})^*$. Since $\mathbf{D}^*$ is locally trivial, there exists a graded vector space $\mathbf{W}=\oplus_{i=1}^{n} W_i$ such that any point in $M$ admits an open neighborhood $U$ with $\mathbf{D}^*|_U= U\times \mathbf{W}$, so that 
\begin{equation}\label{eq:localcoord}
        C_\cM |_U= (\Gamma_{\sym\mathbf{E}^*}/\cI_\mu)|_U\cong \Gamma_{(U\times \sym \mathbf{W})} = C_U^\infty\otimes \sym \mathbf{W}.
        \end{equation}
Since morphisms of graded manifolds are morphisms of ringed spaces it is clear that $\cF$ is well defined at the level of morphisms. 
\end{proof}

\begin{remark} \label{rem:properties}

\begin{itemize}
\item[(a)] The functor $\cF$ takes split $n$-coalgebra bundles (Example \ref{splitexample}) to split $n$-manifolds (Example \ref{ex:splitnman}).
\item[(b)]
The correspondence $(\mathbf{E}\to M,\mu) \mapsto \cF(\mathbf{E}\to M,\mu)=\cM$ respects truncations (see \eqref{tower} and Example~\ref{extrunc}), i.e. $$\cF(\mathbf{E}^{\leq k}\to M, \mu^{\leq k})=\cM_k\quad \text{for all}\quad k\in\{0,\cdots, n\}.$$

\end{itemize}
\end{remark}

\begin{remark}[Local coordinates]\label{rem:local}
For an $n$-manifold $\cM$ as in \eqref{eq:M}, by locally splitting $E^*_{-i}$ as $(\ker \mu_{-i})^* \oplus (K^\mu_{-i})^*$,
we see  that degree $i$ coordinates on $\cM$ are given by a local frame of $(\ker \mu_{-i})^*$ while the spaces $(K^\mu_{-i})^*= (\im \mu_{-i})^*$ codify coordinates in degree $i$ which are products of coordinates of lower degrees. In fact, note that, globally,
\begin{equation}\label{eq:KC}
C^{i}_{\cM_{i-1}}=\Gamma_{(K^\mu_{-i})^*}, \qquad i=2,\cdots, n.
\end{equation}
One can check this equality by noticing that it holds in the split case (see \eqref{eq:Kmu}) and using Prop.~\ref{char-admis}. \hfill $\diamond$
\end{remark}

\subsection{The equivalence of categories}
We will now show that the functor in Prop.~\ref{F-well-def} is an equivalence of categories.

Given an $n$-manifold $\cM=(M, C_\cM)$ recall from $\S$ \ref{S2.2} that,
for each $i=1,2,\ldots$, $C^i_\cM\cong\Gamma_{E_{-i}^*}$ for some vector bundle $E_{-i}\to M$ (unique up to isomorphism), and the multiplication on $C_\cM$ is codified by vector bundle maps $E_{-(i+j)}\to E_{-i}\otimes E_{-j}$.
Picking those maps with $i+j\leq n$ gives a coalgebra structure $\mu$ on $\mathbf{E}:=\oplus_{i=-n}^{-1} E_{-i}$. Note that, by its very construction, this coalgebra bundle $(\mathbf{E},\mu)$ associated with $\cM$ is uniquely defined, up to isomorphism.

\begin{example}\label{ex:split-split}
If $\cM=(M, C_\cM)$ is a split $n$-manifold, then a corresponding coalgebra bundle $(\mathbf{E},\mu)$ is isomorphic to a split coalgebra bundle. Indeed (see Example \ref{ex:splitnman}), since 
$$
C_\cM=\Gamma_{\sym({\bf W})}
$$ 
for a graded vector bundle ${\bf W}=\oplus_{i=1}^nW_i\to M$, by setting ${\bf D}= {\bf W}^*$
one can directly check that  $(\mathbf{E},\mu)\cong (\mathbf{E_D},\mu^s)$ (see Example \ref{splitexample}). \hfill $\diamond$
\end{example}

\begin{lemma}\label{lem:reconst}
Let $\cM$ be an $n$-manifold, and let $(\mathbf{E},\mu)$ be a corresponding coalgebra bundle. Then  $(\mathbf{E},\mu)$ is admissible and $C_\cM\cong \Gamma_{\sym \mathbf{E}^*}/\cI_\mu$.
\end{lemma}

\begin{proof}
Since admissibility is a pointwise condition, it is enough to verify that $(\mathbf{E},\mu)$ is admissible on small open subsets, and this holds since $\cM$ is locally isomorphic to a split $n$-manifold (by \eqref{locality}), so $(\mathbf{E},\mu)$ is locally isomorphic to a split coalgebra bundle (Example~\ref{ex:split-split}), hence admissible (see Example \ref{splitexample}).

For the second statement, note that, by the way $(\mathbf{E},\mu)$ is defined, we have that $C^i_{\cM}=\Gamma_{E_{-i}^*}$, for $i= 1,\cdots, n$, and there is a natural morphism of graded algebras
$\phi:\Gamma_{\sym \mathbf{E}^*}\to C_\cM$ 
satisfying $\phi(\cI_\mu)=0$.
To show that the induced map 
$$
\Gamma_{\sym \mathbf{E}^*}/\cI_\mu \to C_\cM
$$
is an isomorphism, it is enough to show that it is an isomorphism on small open subsets. This can be verified again using that  $(\mathbf{E},\mu)$ is locally isomorphic to a split coalgebra bundle, so the result follows from   Lemma \ref{lemma2}.
\end{proof}

We are ready to show that, for each $n$, the category of $n$-coalgebra bundles and the category of $n$-manifolds are equivalent.

\begin{theorem}\label{geometrization}
For each non-negative integer $n$ the functor $\cF: \cob{n}\to\man{n}$ from Proposition~\ref{F-well-def} is an equivalence of categories. 
\end{theorem}

\begin{proof}
We must verify that $\cF$ is fully faithful and essentially surjective. 
To see that it is fully faithful, let
$(\mathbf{E}\to M,\mu)$ and $(\mathbf{F}\to N, \nu)$ be admissible $n$-coalgebra bundles, and let $\cM=(M,C_\cM)$ and $\cN=(N,C_\cN)$ be their respective images under $\cF$. As noticed in Remark \ref{lowdeg},
we have that $C_\cM^i=\Gamma_{E^*_{-i}}$ and $C_{\cN}^i=\Gamma_{F^*_{-i}}$ for $1\leq i\leq n$. For morphisms $\Phi_1,\Phi_2:(\mathbf{E},\mu)\to (\mathbf{F},\nu)$, consider $\cF(\Phi_1)=(\phi_1,\phi^\sharp_1),\cF(\Phi_2)=(\phi_2,\phi^\sharp_2):\cM\to \cN$, the corresponding morphisms of $n$-manifolds. Then $\cF(\Phi_1)=\cF(\Phi_2)$ if and only if
\begin{equation*}
\left\{\begin{array}{l}
\phi_1=\phi_2,\\
(\phi^{\sharp}_1)^i=(\phi_2^{\sharp})^i:C_\cN^i\to \phi_{1*} C^i_\cM\quad \text{for } i= 1,\ldots, n.
\end{array}\right.
\end{equation*}
 But $(\phi^\sharp_j)^i=(\Phi^*_j)^i:\Gamma_{F_i^*}\to\phi_{j*}\Gamma_{E_i^*}$, for $j=1,2$. Therefore  $\Phi_1=\Phi_2$, so $\cF$ is faithful.

To see that $\cF$ is full, recall that, for a morphism $\Psi=(\psi,\psi^\sharp):\cM\to\cN$, the fact that $\psi^\sharp:C_\cN\to \psi_*C_\cM$ is a morphism of sheaves of graded algebras implies that, for each $i=1,\ldots,n$, $(\psi^\sharp)^i:C^i_\cN=\Gamma_{F_i^*}\to \psi_*C_\cM^i=\psi_*\Gamma_{E_i^*}$  is a morphism of sheaves of $C^\infty_N$-modules, which is equivalent (by adjointness) to a morphism $(\psi^\sharp)^i: \psi^* \Gamma_{F_i^*} \to \Gamma_{E_i^*}$ of sheaves of $C^\infty_M$-modules, which is necessarily defined by a vector-bundle map $\psi^* F_i^*\to E_i^*$. It follows that there is a map $\Phi:\mathbf{E}\to \mathbf{F}$ of graded vector bundles covering $\psi:M\to N$ such that $(\Phi^*)^i=(\psi^\sharp)^i$, for $1\leq i\leq n$. Finally the fact that $\psi^\sharp$ is an algebra morphism implies that $\Phi:\mathbf{E}\to \mathbf{F}$ is a coalgebra morphism.

The fact that $\cF$ is essentially surjective is a direct consequence of Lemma \ref{lem:reconst}.
\end{proof}

As a direct consequence of the equivalence in Theorem~\ref{geometrization}, along with Remark~\ref{rem:properties} (a) and Prop.~\ref{char-admis}, we have the following result, see \cite[Thm.~1]{bon:on}.

\begin{corollary}
Any $n$-manifold is isomorphic to one that is split.    
\end{corollary}

\begin{remark}[Submanifolds]
By means of the geometrization functor, one obtains a description of submanifolds of a given $n$-manifold $\cM=\cF(\mathbf{E}\to M,\mu)$ as admissible coalgebra subbundles $(\mathbf{F}\to N, \nu) \hookrightarrow (\mathbf{E}\to M,\mu)$. Alternatively, by taking annihilators one can describe submanifolds of $\cM$ as graded subbundles $(\mathbf{K}\to N)\subseteq (\mathbf{E}^*\to M)$ satisfying 
\begin{equation*}
     \mathbf{K}\cap \mu^*(\mathbf{E}^*\otimes\mathbf{E}^*)_{|N}=\mu^*(\mathbf{K}\otimes\mathbf{E}^*_{|N}).
    \end{equation*}
At the graded geometric level, such subbundle $\mathbf{K}\to N$ codifies the subsheaf of regular ideals in $C_\cM$ defining the submanifold; see \cite[$\S$ 2.3]{bur:super} for a discussion for the case $n=2$. 
\hfill $\diamond$
\end{remark}

\subsection{Examples}\label{sec:1-tan-cotan}

\begin{example}\label{geo1}
As we saw in Example \ref{cob1}, $\cob{1}$ is isomorphic to $\mathcal{V}ect$, the category of vector bundles.  Theorem \ref{geometrization} boils down to the well-known fact that 
$$
\cF:\mathcal{V}ect\to \man{1}, \quad \cF(E\to M)=(M, \Gamma_{\wedge^\bullet E^*})=E[1]
$$ is an equivalence of categories. \hfill $\diamond$
\end{example}

\begin{example}[$1$-manifolds as $n$-manifolds]\label{vectcoal}
For a vector bundle $E\to M$ define the  admissible $n$-coalgebra bundle $(\mathbf{E}\to M, \mu^E)$ given by $E_{-i}=\wedge^i E$ for $1\leq i\leq n$ and $\mu^E$ the coalgebra structure  dual to the algebra structure on $\mathbf{E}^*$ defined by the wedge product. Then $\cF({\bf E}\to M, \mu^E)=E[1]$, where $E[1]$ is viewed as an $n$-manifold. 
\hfill $\diamond$
\end{example}

\begin{example}
We saw in Example \ref{cob2} a characterization of admissible 2-coalgebra bundles  as triples $(E_{-1}\to M, E_{-2}\to M, \mu:E_{-2}\twoheadrightarrow E_{-1}\wedge E_{-1})$. The geometrization functor  
$$\cF:\cob{2}\to \man{2}, \quad \cF(\mathbf{E}\to M, \mu)=\Big(M, \ \frac{\Gamma_{\wedge^\bullet E_1^*\otimes \sym E_2^*}}{\langle e\wedge e'\otimes 1-1\otimes \mu^*(e,e')\rangle}\Big)
$$
recovers the vector-bundle description of $2$-manifolds used in \cite[$\S$ 2.2]{bur:super}. 
\hfill $\diamond$
\end{example}

The previous examples provided geometric characterizations of $n$-manifolds in low degrees. We now give a class of examples of admissible $n$-coalgebra bundles in arbitrary degrees and describe the 
corresponding $n$-manifolds via the geometrization functor.

In order to construct $n$-coalgebra bundles, we will take as input exact sequences of the form
\begin{equation}\label{eq:dvbs}
0\to C\to \Omega \stackrel{\varphi}{\to} A\otimes B \to 0,
\end{equation}
where $A$, $B$, $C$ and $\Omega$ are vector bundles over $M$, and the maps are over the identity map. We will refer to such exact sequences as {\em dbv-sequences}, since, as shown in \cite{che:ond}, the natural category they form is equivalent to the category of {\em double vector bundles}; more concretely, a double vector bundle
\begin{equation}\label{dbv}
    \xymatrix{D\ar[d]\ar[r]& B\ar[d]\\ A\ar[r]& M}
\end{equation}
with core $C$  gives rise to a dvb-sequence \eqref{eq:dvbs} where $\Omega\to M$ is the dual of the vector bundle whose space of sections are the double linear functions on $D$ (see  \cite{gra:lie} and \cite{mei:wei}), and this sequence fully recovers $D$, up to isomorphism (see also \cite{fer:geo} for a detailed discussion).

For the double vector bundles given by tangent and cotangent prolongations of a given vector bundle $A\to M$, 
\begin{equation}\label{eq:tgcotg}
    \begin{array}{c}
         \xymatrix{TA\ar[d]\ar[r]& TM\ar[d]\\ A\ar[r]& M}
    \end{array}\qquad \text{and}\qquad\begin{array}{c}
         \xymatrix{T^*A\ar[d]\ar[r]& A^*\ar[d]\\ A\ar[r]& M,}
    \end{array} 
\end{equation}
the corresponding dvb-sequences are the dual sequences of
\begin{equation}\label{eq:j1der}
    0\to \Hom(TM,A^*) \to J^1A^*\to A^*\to 0  \qquad\text{and}\qquad 0\to \End(A) \to \der(A)\to TM\to 0,
\end{equation}
where $J^1A^*$ is the \emph{first jet prolongation} of $A^*$ and $\der(A)$ is the vector bundle whose sections are derivations of $A$ (see $\S$ \ref{G-vf}).

\begin{example}[$n$-coalgebra bundles from double vector bundles] \label{ex:dvb}
For $D$ a double vector bundle (as in \eqref{dbv}) with corresponding dvb-sequence \eqref{eq:dvbs}, we define, for each $n\geq 2$,
an admissible $n$-coalgebra bundle $(\mathbf{E}^D\to M, \mu^D)$ as follows. 
For $n=2$, we set
$$
E^D_{-1}=A\oplus B,\quad E^D_{-2}=\wedge^2 A\oplus \Omega \oplus \wedge^2 B,\quad \mu^D_{-2}(\xi+ \omega+\tau)=\mu^A_{-2}(\xi)+\varphi(\omega)+\mu_{-2}^B(\tau),
$$
where $\xi\in\wedge^2A,\ \omega\in \Omega,\ \tau\in\wedge^2 B$ and $\mu^A,\mu^B$ are as in Example \ref{vectcoal}. For $n>2$, we set 
$$
E^D_{-i}=\left\{\begin{array}{ll}
    \wedge^iA, &  i=1,\ldots, n-2,\\
    \wedge^{n-1} A\oplus B, & i=n-1,\\
    \wedge^n A\oplus \Omega, & i=n,
\end{array}\right.\quad\text{and}\quad \mu^D_{-i}=\left\{\begin{array}{ll}
   \mu^A_{-i}, &  i=1,\ldots, n-1\\
    \mu^A_{-n}+\varphi, & i=n.
\end{array}\right.
$$
Note that one can directly check that $(\mathbf{E}^D\to M, \mu^D)$ is admissible by  realizing that an splitting of the exact sequence \eqref{eq:dvbs} gives an isomorphism with the split coalgebra bundle  $(\mathbf{E_D}\to M, \mu^s)$, where $\mathbf{D}\to M$ is given by $D_{-1}=A$, $D_{1-n}=B$ and $D_{-n}=C$. 

For a geometric description of the $n$-manifolds 
$$
\cM^D:= \cF(\mathbf{E}^D\to M, \mu^D),
$$
we consider the $1$-manifolds $A[1]$ and $D[1]$ corresponding to the vector bundles $A\to M$ and $D\to B$. The additional vector bundle structure of $D$ over $A$ makes $D[1]\to A[1]$ into a vector bundle in the category of $1$-manifolds. Denoting by $\mathcal{R}\to A[1]$ the trivial vector bundle with just one fiber coordinate of degree $n-1$, i.e. $\mathcal{R}=\bR^{0|\cdots|0|d_{n-1}}\times A[1]$ with $d_{n-1}=1$, the $n$-manifold $\cM^D$ can be regarded as the total space of the graded vector bundle $D[1]\otimes \mathcal{R}\to A[1]$. 

When $D$ is the tangent or cotangent prolongation of $A\to M$ (see \eqref{eq:tgcotg} and \eqref{eq:j1der}), one can check that  
$$
\cM^{TA}=T[n-1]A[1] \quad \text{and}\quad \cM^{T^*A}=T^*[n]A[1],
$$ 
as defined in Example \ref{shift-tan-cotan}.
\hfill $\diamond$
\end{example}

One can also use coalgebra bundles to give more explicit descriptions of morphisms between $\mathbb{N}$-manifolds, see e.g. Remark~\ref{rem:posdegvf}.

\section{Vector fields and tangent vectors}\label{S4}

In this section we review the notion of vector fields on $\mathbb{N}$-manifolds and characterize them in terms of admissible coalgebra bundles.

\subsection{Vector fields}
Let $\cM=(M,C_\cM)$ be an $n$-manifold of dimension $m_0|\cdots|m_n$. For an open subset $U\subseteq M$, a \emph{vector field of degree $k$} on $\cM|_U$ is a degree $k$ derivation $X$ of $C_\cM(U)$, i.e.,  an  $\mathbb R$-linear map $X:C_\cM(U)\to C_\cM(U)$ with the property that, for all $f, g \in C_\cM(U)$ with $f$ homogeneous, $|X(f)| =|f|+k$ and
    \begin{equation}\label{der}
 X(fg)=X(f)g+(-1)^{|f|k}fX(g).
    \end{equation}
Vector fields give rise to a sheaf of $C_\cM$-modules over $M$. 
The sheaf of degree $k$ vector fields is denoted by $\vf^k$, and the sheaf of all vector fields by $\vf^\bullet= \bigoplus_{k \in \mathbb{Z}} \vf^k$. The graded commutator of vector fields, defined for homogeneous vector fields $X, Y$ by
\begin{equation*}
    [X,Y]=XY-(-1)^{|X||Y|}YX,
\end{equation*}
makes $\vf^{\bullet}$ into a sheaf of graded Lie algebras.

Let $U\subseteq M$ be a chart of $\cM$ with coordinates $\{ x^\alpha, e^{\beta_i}_i\}$. Define the vector fields $\frac{\partial}{\partial x^\alpha},\ \frac{\partial}{\partial e^{\beta_i}_i}$ as the derivations acting on coordinates by
$$\frac{\partial}{\partial x^\alpha}(x^{\alpha'})=\delta_{\alpha \alpha'}, \quad \frac{\partial}{\partial x^\alpha}(e_i^{\beta_i})=0,\quad \frac{\partial}{\partial e^{\beta_i}_i}(x^\alpha)=0,\quad \frac{\partial}{\partial e^{\beta_i}_i}(e^{\beta'_j}_j)=\delta_{ij}\delta_{\beta_i \beta'_j}.$$
This  definition implies that $\big|\frac{\partial}{\partial x^{\alpha}}\big|=0$ and $\big |\frac{\partial}{\partial e^{\beta_i}_i}\big |=-\big |e^{\beta_i}_i \big |=-i$.

\begin{proposition}\label{neg-gen-vf}
Let $U$ be a chart of $\cM$ with coordinates $\{ x^\alpha, e^{\beta_i}_i\}$. Then
$\vf^{\bullet}(U)$ is freely generated by $\{\frac{\partial}{\partial x^{\alpha}}, \frac{\partial}{\partial e^{\beta_i}_i}\}$ as a $C_{\cM}(U)\text{-module}$.
\end{proposition}
\begin{proof} 
Let  $X$ be a vector field on $\cM_{|U}$, and consider the vector field
\begin{equation*}
X'=X-\sum_{\alpha=0}^{m_0} X(x^\alpha)\frac{\partial}{\partial x^{\alpha}}- \sum_{i=1}^{n}\sum_{\beta_i=1}^{m_i} X(e^{\beta_i}_i)\frac{\partial}{\partial e^{\beta_i}_i}.
\end{equation*}
It is clear that $X'$ acts by zero on the coordinates, hence on arbitrary functions, and therefore $X$ can be written as a $C_\cM(U)$-linear combination of $\{\frac{\partial}{\partial x^{\alpha}}, \frac{\partial}{\partial e^{\beta_i}_i}\}$.
\end{proof}

\begin{remark}\label{vf-sec-tan}
By considering vector bundles in the realm of graded manifolds,
one can equivalently view vector fields on $\cM$ 
as sections of tangent bundles \cite{raj:tes}.
\hfill $\diamond$
\end{remark}

\begin{example}[$Q$-manifolds]\label{ex:Qman}
A degree $1$ vector field $Q$ on an $n$-manifold $\cM$ such that $[Q,Q]=2Q^2=0$ is called a \emph{homological vector field} (aka a {\em $Q$-structure}), since
in this case $(C_\cM, Q)$ is a differential complex. The pair $(\cM , Q)$ is called a \emph{$Q$-manifold}. A 1-manifold with a homological vector field is the same as a Lie algebroid \cite{vai:lie}; in general, $Q$-manifolds codify $L_n$-algebroids, see \cite{bon:on}.
\hfill $\diamond$
\end{example}

\subsection{Tangent vectors at a point}
Let $\cM=(M,C_\cM)$ be an $n$-manifold of dimension $m_0|\cdots|m_n$. 
Given $x\in M$, we denote by $C_{\cM}|_x$ the \emph{stalk} of $C_\cM$ at $x$, and given $f\in C_\cM(U)$ with $x\in U$, we denote by $\mathbf{f}$ its class in $C_{\cM}|_x$. 
A \emph{homogeneous tangent vector of degree $k$ at a point $x\in M$} is 
a linear map $v:C_{\cM}|_x\to \bR$ satisfying
\begin{equation*}
 v(\mathbf{fg})=v(\mathbf{f})g^0(x)+(-1)^{k|f|}f^0(x)v(\mathbf{g}),
\end{equation*}
where $f^0$ and $g^0$ denote the degree zero components of $f$ and $g$, respectively. We denote by $T_x\cM$ the space of tangent vectors at a point $x\in M$, noticing that it is a graded vector space over $\mathbb{R}$ with $\dim (T_x\cM)_{-i}=m_i$ for $i=0,\ldots, n$, and zero otherwise.

For each open subset $U\subseteq M$, any vector field $X\in\vf^{\bullet}(U)$ defines a tangent vector $X_x$ at each point $x\in U$ by
\begin{equation}\label{tan-vf-gen}
X_x(\mathbf{f})=X(f)^0(x).
\end{equation}

\begin{definition}\label{def:li} We say that two vector fields $X,Y\in\vf^{\bullet}(U)$ are \emph{linearly independent} if $X_x$ and $Y_x$ are linearly independent for all $x \in U$. 
\end{definition}

For $\mathbb{N}$-manifolds, it is clear that only  vector fields of non-positive degrees can generate non-zero tangent vectors. Therefore in this setting, unlike the case of smooth manifolds, vector fields are not determined by their corresponding tangent vectors.

\begin{example}[Two different vector fields with the same tangent vectors at all points]\label{Vf-tang}
Consider the $1$-manifold $\bR^{1|1}:=\{ x, e\}$, where $|x|=0$ and $|e|=1$. Then the vector fields 
$$
X=\frac{\partial}{\partial x} \quad\text{and}\quad Y=\frac{\partial}{\partial x}+e \frac{\partial}{\partial e}
$$ 
have the same tangent vectors at all points. \hfill $\diamond$
\end{example}

\begin{proposition}\label{prop:li}
Let $U\subseteq M$ be an open subset, and consider vector fields $X_k^{i_k}\in \cT_{\cM}^{-k}(U)$, for $k=0,\ldots, n$, $i_k=1,\ldots, d_k$.
If they are linearly independent (in the sense of Definition~\ref{def:li}), then the set $\{X_k^{i_k}\}$ is linearly independent over $C_\cM(U)$.
\end{proposition}

The proof is a direct extension of the argument in \cite[Prop.~2.17]{bur:super} for $n=2$.

\subsection{Geometrization of vector fields}\label{G-vf}
We will now describe vector fields in non-positive degrees on $\mathbb{N}$-manifolds in terms of coalgebra bundles.

Let $\mathbf{E}=\oplus_{i=k}^l E_i\to M$ be a graded vector bundle, for $k,l\in\bZ$. For an open subset $U\subseteq M$, a {\em derivation} of $\mathbf{E}|_U$ is a pair 
$(D,X)$ where $X$ is a vector field on $U$, called the {\em symbol} of the derivation, and $D: \Gamma_{\mathbf{E}}(U) \to \Gamma_{\mathbf{E}}(U)$ 
satisfies
\begin{equation*}
D(fe)=X(f)e+fD(e), \qquad \; \; \forall \; f\in C^\infty(U), \; e\in\Gamma_\mathbf{E}(U).
\end{equation*}
Note that if $(D, X)$ is  homogeneous of nonzero degree, then $X = 0$. The corresponding sheaf of derivations of $\mathbf{E}\to M$ is a sheaf of (graded) $C^\infty_M$-modules and also a sheaf of graded Lie algebras with respect to 
the graded commutator bracket
\begin{equation*}
[D,D']=D\circ D'-(-1)^{|D||D'|} D'\circ D.
\end{equation*}
Since, as a sheaf of $C^\infty_M$-modules, the sheaf of derivations is locally finitely generated and free, it can be realized as the sheaf of sections of a (graded) vector bundle 
$$
\der(\mathbf{E})\to M
$$ 
fitting into the following short exact sequence,
\begin{equation}\label{eq:atiyah}
0\to \End(\mathbf{E}) \to \der(\mathbf{E}) \stackrel{\sigma}{\to} TM \to 0,
\end{equation}
where $\sigma$ is the {\em symbol map}.
Note that $\der(\mathbf{E})=\oplus_j \der(\mathbf{E})_j$, where 
\begin{equation}\label{atiyah}
\der(\mathbf{E})_j=\begin{cases}
\der(E_k)\times_{TM}\cdots\times_{TM}\der(E_l)& \text{for } j=0,\\
\bigoplus\limits_{b=a+j} \Hom(E_a, E_b)& \text{for }j\neq 0.
\end{cases}
\end{equation}

\begin{remark}
The bundle $\der(\mathbf{E})$, also called the {\em Atiyah algebroid} of $\mathbf{E}\to M$, was previously considered in \cite{raj:reprep}, where representations up to homotopy of a Lie algebroid $\Alie$ on a graded vector bundle $\mathbf{E}\to M$ are understood as  $L_\infty$-morphisms from $A$ to $\der(\mathbf{E})$. \hfill $\diamond$
\end{remark}

We will use such derivations to describe vector fields of {non-positive} degrees on $\mathbb{N}$-manifolds. 
Given an $n$-coalgebra bundle $(\mathbf{E}\to M, \mu)$, we  define the associated  \emph{extended $n$-algebra bundle} as the pair $(\mathbf{E}_{\bR}^*\to M, m)$, where 
$\mathbf{E}_{\bR}^* = \oplus_{i=0}^n (\mathbf{E}_\bR^*)_i$ is the graded vector bundle with
\begin{equation}\label{extended}
(\mathbf{E}_\bR^*)_i=\left\{\begin{array}{ll}
\bR_M,& i=0,\\
E_{-i}^*,& 1\leq i\leq n,
\end{array}\right. 
\end{equation}
and the multiplication $m: \mathbf{E}_{\bR}^*\otimes \mathbf{E}_{\bR}^*\to \mathbf{E}_{\bR}^*$ is defined by $\mu^*$ and scalar multiplication $\bR_M\times E_j^*\to E_j^*$; here $\mathbb{R}_M=\mathbb{R}\times M \to M$. For an open subset $U\subseteq M$, a degree $k$ derivation $(D, X)$ of $\mathbf{E}_\bR^*|_U$ is said to be {\em compatible} with $m$ if 
\begin{equation}\label{der-m}
D(m(e,e'))=m(D(e),e')+(-1)^{k|e|}m(e,D(e')).
\end{equation}
One can directly verify  that the sheaf of compatible derivations of $(\mathbf{E}_{\bR}^*\to M, m)$ is a subsheaf of graded Lie algebras of the sheaf of derivations of $\mathbf{E}_{\bR}^*$.

\begin{lemma}\label{geo-vf}
Let $(\mathbf{E}\to M,\mu)$ be an admissible $n$-coalgebra bundle and $\cM=\cF(\mathbf{E}\to M,\mu)$ be the corresponding $n$-manifold. For each non-positive integer $k$,  there is a natural inclusion 
$$
\Upsilon^k: \cT^k_\cM  \hookrightarrow \Gamma_{\der(\mathbf{E}_\bR^*)_k},
$$ 
given, for each open subset $U\subseteq M$, by
\begin{equation*}
 \Upsilon^k(X)(e)=X(e),\qquad \Upsilon^k(X)(f)=X(f),
\end{equation*}
where $X\in \cT^k_\cM(U)$,  $e\in\Gamma_{(\mathbf{E}^*_{\bR})_i}(U)=\Gamma_{E^*_{-i}}(U)=C^i_\cM(U)$, for $i=1,\ldots,n$, and  $f\in\Gamma_{\bR_M}(U)=C^\infty(U)$, with the following properties:
 \begin{enumerate}
     \item 
     $\Upsilon^k(\cT^k_\cM(U)) =\{D\in\Gamma_{\der(\mathbf{E}_\bR^*)_k}(U) \text{ compatible with $m$} \}, $
     \item $\Upsilon^{k+l}([X, Y])=[\Upsilon^k(X), \Upsilon^l(Y)]$.
 \end{enumerate}

\end{lemma}

\begin{proof}
The $n$-manifold $\cM =  (M, C_\cM)$ is such that $C^i_\cM=\Gamma_{(\mathbf{E}^*_{\bR})_i}$ for $i=0,\ldots,n$, and the multiplication of functions in $C_\cM^i$ and $C_\cM^j$, with $i+j\leq n$,  agrees with $m$. Since $C_\cM$ is generated as a sheaf of algebras by the homogeneous functions of degree $\leq n$, for any open subset $U\subseteq M$, a vector field  $X\in\vf^{k}(U)$ is completely determined by the maps
\begin{equation*}
 X^i:C^i_\cM(U)\to C^{i+k}_\cM(U) \quad \text{for } i=0,\ldots,n.
\end{equation*}
Thus, for $k\leq 0$, we see that $X$ gives rise, and is completely determined by, a map $D:\Gamma_{\mathbf{E}^*_{\bR}}(U)\to \Gamma_{\mathbf{E}^*_{\bR}}(U)$, which is a degree $k$ derivation. Hence $\Upsilon^k$ is well defined and injective.  The bracket preserving property (b) follows directly from the definitions. 

It remains to verify (a). Since $m$ agrees with the multiplication of functions on $C_\cM$, it is clear that any derivation in the image of $\Upsilon^k$ satisfies \eqref{der-m}, so we obtain one inclusion in (b). For the other inclusion observe that any degree $k$ derivation $D$ of $\mathbf{E}^*_{\bR}$ can be extended to $\widehat{D}:\Gamma_{\sym \mathbf{E}^*}\to \Gamma_{\sym \mathbf{E}^*}$ as a degree $k$ derivation. If $D$ satisfies \eqref{der-m} then $\widehat{D}(\cI_\mu)\subseteq \cI_\mu$, so it descends to a derivation $X$ of $C_\cM^\bullet=\Gamma_{\sym \mathbf{E}^*}/\cI_\mu$, in such a way that $\Upsilon^k(X)=D$. 
\end{proof}

Since $\vf^k$ is locally finitely generated and free as a sheaf of $C^\infty_M$-modules, a consequence of Theorem \ref{geo-vf} is that the sheaf of compatible derivations of $(\mathbf{E}^*_{\bR}, m)$ of degree $k$ is the sheaf of
sections of a subbundle 
$$
\der(\mathbf{E}^*_{\bR},m)_k \subseteq \der(\mathbf{E}^*_{\bR}),
$$
for $k=-n,\ldots, 0$. This subbundle has a simple description when $k=-n$.
From \eqref{atiyah}, note that $\der(\mathbf{E}^*_{\bR})_{-n}=\Hom(E_{-n}^*,\bR)=E_{-n}$, while from \eqref{der-m} we see that $\xi\in E_{-n}$ is a compatible derivation if and only if $\xi|_{(K^\mu_{-n})^*}=0$ (see \eqref{eq:KC}), which implies that
$$
\der(\mathbf{E}^*_{\bR},m)_{-n} = \ker \mu_{-n} \subseteq E_{-n}.
$$

\begin{theorem}\label{thm:compatder}
For each $k=-n,\ldots, 0$, 
\begin{equation*}
\vf^k = \Gamma_{\der(\mathbf{E}^*_{\bR},m)_k}.
\end{equation*}
In particular, $\vf^{-n}=\Gamma_{\ker \mu_{-n}}$.
\end{theorem}

By means of the identification in Theorem~\ref{thm:compatder}, the natural $C_\cM$-module structure of $\cT_\cM$ translates into  vector bundle maps
\begin{equation}\label{eq:theta}
\Theta_{i,k}: E_{-i}^*\otimes \der({\bf E}^*_\bR,m)_{k} \to \der({\bf E}^*_\bR,m)_{k+i}, 
\end{equation}
given on sections by $\Theta_{i,k}(e\otimes D)=m(e, D(\cdot))$, for $i=1,\ldots,n$, $k=-n,\ldots,0$, and $k+i\leq 0$.

\begin{example}[Vector fields on $1$-manifolds]\label{1-vf}
Let $\cM=E[1]$ be a $1$-manifold. As recalled in Example \ref{cob1}, the associated coalgebra bundle is $\mathbf{E} = E_{-1}=E$, with zero comultiplication. By  \eqref{atiyah} and Theorem \ref{thm:compatder},
$$
\der(\mathbf{E}^*_{\bR},m)_{0}=  {\der(E^*)}\cong {\der(E)}, \qquad 
\der(\mathbf{E}^*_{\bR})_{-1}=E,
$$
where the isomorphism on the left-hand side is via dualization of derivations.
Therefore
\begin{equation*}
    \cT^0_{E[1]} = \Gamma_{\der(E^*)}\cong \Gamma_{\der(E)} 
    \quad \text{ and } \quad \cT^{-1}_{E[1]} = \Gamma_{E},
\end{equation*}
see \cite{zam:dis}. The only nontrivial map in the collection \eqref{eq:theta} occurs when $i=1$ and $k=-1$, and it is the natural inclusion
$$
E^*\otimes E=\End(E)\hookrightarrow \der(E),
$$
see \eqref{eq:atiyah}.
\hfill $\diamond$
\end{example}

We will discuss vector fields on 2-manifolds in Example~\ref{2-vf} below.

\begin{remark}[Positive-degree vector fields on $1$-manifolds]\label{rem:posdegvf}
One can also describe  positive-degree vector fields on $E[1]$. For $k>0$, since a degree $k$ vector field  $Q$ is a degree $k$ derivation of $\Gamma_{\wedge^\bullet E^*}$, it is completely described by maps 
$$
\delta_0: C^\infty_M \to \Gamma_{\wedge^k E^*} \;\;  \mbox{ and } \;\;  \delta_1: \Gamma_{E^*}\to \Gamma_{\wedge^{k+1} E^*}
$$
satisfying $\delta_{0}(fg)=\delta_0(f)g + f\delta_0(g)$ and $\delta_1(f e)= \delta_0(f)\wedge e + f \delta_1(f)$. Such pairs $(\delta_0,\delta_1)$ are called {\em almost $(k+1)$-differentials} in \cite{xu:uni}, where they are shown to be equivalent to linear $(k+1)$-vector fields on $E^*$.  For $k=1$ and $Q^2=0$, this recovers the correspondence between $Q$-structures on $E[1]$, Lie-algebroid structures on $E$  (see Example~\ref{ex:Qman}) and linear Poisson structures on $E^*$.

From another viewpoint,  following Remark \ref{vf-sec-tan}, a degree $k$ vector field $Q$, with $k>0$, may be seen as a section 
$Q:E[1]\to T[k]E[1]$
of the shifted tangent bundle of $E[1]$ (here $Q$ is regarded a morphism of $(k+1)$-manifolds). Using the geometrization functor in Examples~\ref{vectcoal} and \ref{ex:dvb}, one can view $Q$ as a morphism $Q:({\bf E},\mu^E)\to({\bf E}^{TE},\mu^{TE})$ of $(k+1)$-coalgebra bundles satisfying $\id_{\bf E}=\pi\circ Q$, where $\pi:({\bf E}^{TE},\mu^{TE})\to ({\bf E},\mu^E)$ is the natural projection. We conclude that, for $k\geq 1$, there is a natural identification 
\begin{equation*}
  \mbox{degree $k$ vector field }\,  Q \rightleftharpoons\left\{\begin{array}{l}
        \text{Pairs of vector bundle maps}\\
         \rho:\wedge^{k}E\to TM, \text{ and }
         \Theta:\wedge^{k+1}E\to (J^1E^*)^*\\
         \text{ such that }\varphi\circ\Theta=\rho\wedge\id
    \end{array}\right.
\end{equation*}
By means of the canonical splitting $\Gamma_{J^1E^*}=\Gamma _{E^*}\oplus (\Omega^1_M\otimes \Gamma_{E^*})$, the previous description of vector fields in positive degrees can be alternatively expressed in terms of multi-brackets,
\begin{equation*}
    \mbox{degree $k$ vector field }\,  Q\rightleftharpoons\left\{\begin{array}{l}
        \text{A vector bundle map } \rho:\wedge^{k}E\to TM,\\
          \text{and a bracket }
         [\cdot,\cdots,\cdot]:\Gamma_{\wedge^{k+1} E}\to \Gamma_E \text{ such that}\\ \text{}
         [e_1,\cdots,fe_{k+1}]=f[e_1,\cdots,e_{k+1}]+\rho(e_1\wedge\cdots\wedge e_k)(f)e_{k+1}.
    \end{array}\right.
\end{equation*}
The fact that the objects on the right-hand side are equivalent to almost differentials on $E^*\to M$ is verified in \cite[Sec.~ 3.2]{xu:uni}.
\hfill $\diamond$
\end{remark}

\subsection{An inductive description}\label{subsec:induc}
Recall (see \eqref{tower}) that an $n$-manifold $\cM$ gives rise to a tower of $\mathbb{N}$-manifolds
$$
M=\cM_0\leftarrow\cM_1\cdots\leftarrow\cM_{n-1}\leftarrow\cM_{n}=\cM,
$$
where $\cM_r$ is an $r$-manifold defined by the subsheaf of algebras of $C_\cM$ locally generated by functions of degree $\leq r$. 
We will use this tower to obtain a more explicit description of vector fields in non-positive degrees on $\cM$ following Theorem \ref{geo-vf}.

For $k\leq 0$, any section of $\cT^k_\cM$ preserves the subsheaf $C_{\cM_r}\subseteq C_\cM$, and hence restricts to a section of $\cT^{k}_{\cM_r}$, so we have restriction maps
$$
\mathrm{res}_r : \cT^k_\cM \to \cT^{k}_{\cM_r},\quad \quad r=0,\cdots, n.
$$

\begin{remark}\label{rem:symbol}
Let $\cM$ correspond to the admissible $n$-coalgebra bundle $(\mathbf{E},\mu)$.
When $k=0$, $r=0$, the restriction map $\cT^0_\cM \to \cT^0_{M}= \mathfrak{X}^1_M$ corresponds to the symbol map 
$$
\der(\mathbf{E}^*_\bR,m)_0\subseteq \der(\mathbf{E}^*_\bR)\stackrel{\sigma}{\to} TM.
$$\hfill $\diamond$
\end{remark}

Let $X$ be a vector field on $\cM$ of degree $k\leq 0$.
Since a vector field on $\cM$ is completely specified by how it acts on $C_{\cM}^i$, for $i=1,\ldots,n$, and 
$$
C_{\cM}^i=C^i_{\cM_{n-1}}, \qquad \mbox{for} \;\; i=1,\ldots, n-1,
$$
it follows that $X$ is completely determined by the degree $k$ vector field $ \mathrm{res}_{n-1}(X)$ on $\cM_{n-1}$ together with the map $\mathrm{p}_n(X):= X_{|C^n_\cM}:C^n_\cM\to C^{n+k}_\cM$. Note that $\mathrm{res}_{n-1}(X)$ and $\mathrm{p}_n(X)$ are related by the fact that they coincide on $C^n_{\cM_{n-1}}\subseteq C_{\cM}^n$. In fact, we have

\begin{lemma}\label{lem:vfs}
A vector field on an $n$-manifold $\cM$ of degree $k\leq 0$ is equivalent to a pair $(Y,Z)$, where $Y$ is a degree $k$ vector field on $\cM_{n-1}$, $Z: C^n_\cM\to C^{n+k}_\cM$ satisfies $Z(fe)=Z(f)e + f Z(e)$ for $f$ and $e$ sections of $C^\infty_M$ and $C^n_\cM$, respectively, and such that 
$$
Z_{|_{C^n_{\cM_{n-1}}}} =  Y_{|_{C^n_{\cM_{n-1}}}}.
$$
\end{lemma}

Let $(\mathbf{E}\to M,\mu)$ be an admissible $(n+1)$-coalgebra bundle, let
$\cM=\cF(\mathbf{E},\mu)$ be the corresponding  $(n+1)$-manifold, take $k\leq 0$, and recall that
$$
\vf^k = \Gamma_{\der(\mathbf{E}^*_{\bR},m)_k}, 
$$
and that
$$
\cT^k_{\cM_n}=\Gamma_{\der((\mathbf{E}^{\leq n})_{\bR}^*, m^{\leq n})_k} \quad \mbox{ and } \quad
C_{\cM_n}^{n+1} =  \Gamma_{(K^\mu_{-(n+1)})^*}, 
$$
see Remarks~\ref{rem:properties} (b) and \ref{rem:local}.
The next result is a reformulation of the previous lemma in terms of coalgebra bundles using the identifications above.

\begin{proposition}\label{prop:pullbacks}
The vector bundle $\der(\mathbf{E}_{\bR}^*, m)_{0}$ is given by the following pullback diagram:
\begin{equation}\label{Geo-deg0-vf}
    \xymatrix{\der(\mathbf{E}_{\bR}^*, m)_{0}\ar[d]^{\mathrm{res}_n}\ar[r]^{\mathrm{p}_{n+1}}& \der(E^*_{-(n+1)}) \ar[d]\\ \der((\mathbf{E}^{\leq n})_{\bR}^*, m^{\leq n})_0\ar[r]& \der((K^\mu_{-(n+1)})^*). }
\end{equation}
Analogously, for $k<0$, the vector bundles $\der(\mathbf{E}_{\bR}^*, m)_{k}$ are given by the pullback diagrams
\begin{equation}\label{Geo-degk-vf}
    \xymatrix{\der(\mathbf{E}_{\bR}^*, m)_{k}\ar[d]^{\mathrm{res}_n}\ar[r]^{\mathrm{p}_{n+1}\qquad}& \Hom(E^*_{-(n+1)}, E^*_{-(n+1)+k}) \ar[d]\\ \der((\mathbf{E}^{\leq n})_{\bR}^*, m^{\leq n})_k\ar[r]&\Hom((K^\mu_{-(n+1)})^* , E^*_{-(n+1)+k}). }
\end{equation}
\end{proposition}

Following the notation above, the maps $\mathrm{res}_n$ and ${\mathrm{p}_{n+1}}$ are given by the natural restrictions. Recalling the exact sequence
$$
0\to \ker \mu_{-(n+1)}\to E_{-(n+1)}\xrightarrow{\mu_{-(n+1)}} \im \mu_{-(n+1)}=K^\mu_{-(n+1)}\to 0,
$$
the right vertical maps in the diagrams are given by restrictions with respect to the inclusions  $\mu^*_{-(n+1)}:(K^\mu_{-(n+1)})^*\hookrightarrow E^*_{-(n+1)}$. 
The bottom maps are defined by  means of \eqref{der-m}.

As an illustration of the previous proposition, we give a description of non-positive degree vector fields on $2$-manifolds using the description of vector fields on $1$-manifolds given in Example \ref{1-vf}.

\begin{example}[Vector fields on $2$-manifolds]\label{2-vf}
Let $(\mathbf{E}=E_{-2}\oplus E_{-1}, \mu: E_{-2}\to E_{-1}\otimes E_{-1})$ be an admissible $2$-coalgebra bundle and denote by $\cM$ the associated degree $2$-manifold. Then putting together Examples \ref{cob2} and \ref{1-vf} we get that $K^{\mu}_{-2}=E_{-1}\wedge E_{-1}$, $\cM_1=E_{-1}[1]$ and $$\cT^{-2}_{\cM_1}=0,\quad \cT^{-1}_{\cM_1}=\Gamma_{E_{-1}}\quad \text{and}\quad \cT^0_{\cM_1}=\Gamma_{\der(E_{-1}^*)}.
$$ 
Thus, for $k=0, -1, -2$, the vector bundles $\der(\mathbf{E}_{\bR}^*, m)_{k}$  are defined by the following pullback diagrams:
\begin{equation*}
    \xymatrix{\der(\mathbf{E}_{\bR}^*, m)_{-2}\ar[d]^{\mathrm{res}_1}\ar[r]^{\quad \mathrm{p}_2}& E_{-2}\ar[d]^{\mu}\\0 \ar[r]^{0\qquad}&E_{-1}\wedge E_{-1,} } \quad 
    \xymatrix{\der(\mathbf{E}_{\bR}^*, m)_{-1}\ar[d]^{\mathrm{res}_1}\ar[r]^{\mathrm{p}_2}& \Hom(E^*_{-2}, E_{-1}^*) \ar[d]^{\mu\otimes \id_{E_{-1}^*}}\\ E_{-1} \ar[r]^{}&\Hom(E_{-1}^*\wedge E_{-1}^* , E_{-1}^*), }
\end{equation*}

\begin{equation*}
    \xymatrix{\der(\mathbf{E}_{\bR}^*, m)_{0}\ar[d]^{\mathrm{res}_1}\ar[r]^{\mathrm{p}_2}& \der(E^*_{-2}) \ar[d]\\ \der(E_{-1}^*) \ar[r]&\der(E_{-1}^*\wedge E_{-1}^*), }
\end{equation*}
with the bottom maps given by 
$$\xi(e\wedge e')=\xi(e) e'- \xi(e') e\quad \text{and} \quad D(e\wedge e')=D(e)\wedge e'+e\wedge D(e'),
$$
for $\xi\in\Gamma({E_{-1}}), \ D\in\Gamma({\der(E_{-1}^*)})$ and $e, e'\in\Gamma(E_{-1}^*)$. In other words,
\begin{align*}
\der(\mathbf{E}_{\bR}^*, m)_{-2} & =\ker (\mu), \\
\der(\mathbf{E}_{\bR}^*, m)_{-1} & = \{ (\xi, \varphi)\in {E_{-1}}\times {\Hom(E^*_{-2}, E_{-1}^*)}\ | \ \varphi(e\wedge e')= \xi(e) e'-\xi(e')  e\},
\end{align*}
and $\der(\mathbf{E}_{\bR}^*, m)_{0}$ is such that
$$
\Gamma(\der(\mathbf{E}_{\bR}^*, m)_{0})=\{ (D^1, D^2)\in \Gamma({\der(E_{-1}^*)\times_{TM}\der(E_{-2}^*)}) \, | \, D^2(e\wedge e')=D^1(e)\wedge e'+e\wedge D^1(e')\}.
$$
This recovers the geometric description of vector fields in degrees $0$, $-1$ and $-2$ in \cite{fer:geo}.
\hfill $\diamond$
\end{example}

\section{Distributions}\label{S5}

Let $\cM=(M,C_\cM)$ be an $n$-manifold of dimension $m_0|\cdots|m_n$. 

\begin{definition}\label{def:dist}
A \emph{distribution} on $\cM$ of rank $d_0|\cdots|d_n$ is a graded\footnote{By graded we mean that if $X$ is a section of $\cD$ over $U$, then so are its homogeneous components.} subsheaf 
$$
\cD \subseteq \vf^{\bullet}
$$
of $C_\cM$-modules such that any point in $ M$ admits an open neighborhood $U$ with the property that $\cD(U)$ is generated by vector fields $\{X_k^{i_k}\}$, where $k=0,\ldots, n$ and $i_k=1,\ldots, d_k$, such that $|X_k^{i_k}|= -k$ and which are linearly independent (at each point in $U$). 

We say that a distribution $\cD$ is \emph{involutive} if it is closed under the Lie bracket on $\vf^{\bullet}$.
\end{definition}

For $n=0$ such distributions are just the usual (regular) distributions on manifolds.

\begin{remark} 
Distributions on supermanifolds have been considered e.g. in 
\cite{fio:sup, cov:fro, del:not} as subsheaves of the tangent sheaf which are locally direct factors. Our formulation of the local triviality property is closer to the one in \cite{lei:sos}.
\end{remark}

\begin{example}[The tangent vectors of a distribution do not determine its involutivity]\label{ex2-dis}
In $\bR^{0|1|1}$, described by coordinates $\{e,p\}$ with $|e|=1$ and $|p|=2$, consider the distributions
\begin{equation*}
\cD=\Big \langle \frac{\partial}{\partial e} \Big \rangle \quad\text{and}\quad \cD'= \Big \langle \frac{\partial}{\partial e}+e\frac{\partial}{\partial p} \Big \rangle.
\end{equation*}
These two distributions have the same tangent vectors at all points, but $\cD$ is involutive while $\cD'$ is not.
\end{example}

\begin{lemma}
For each $l=0,\ldots, n$, the sheaf  of $C^\infty_M$-modules  $\cD\cap \vf^{-l}$ is locally free and finitely generated.
\end{lemma}

\begin{proof}
For a small open subset $U\subseteq M$, 
let $\{X_k^{i_k}\}$ be linearly independent vector fields generating $\cD(U)$ (as in the definition), and, for each $r=1,\ldots, n$, let  $\{ w_r^{q_r}\}$ be a basis of $C_\cM^r(U)$ as a $C^\infty(U)$-module. Then, as a $C^\infty(U)$-module,
\begin{equation}\label{eq:localgen}
\cD\cap \vf^{-l} (U) = \mathrm{span}_{C^\infty(U)}\{ X_l^{i_l}, w_{r}^{q_{r}} X_{l+r}^{i_{l+r}}\, | \, r=1,\ldots, n-l\}.
\end{equation}
The fact that the generating set above is linearly independent over $C^\infty(U)$ is a direct consequence of the linear independence of the vector fields $\{X_k^{i_k}\}$ over $C_\cM(U)$, which follows from Prop.~\ref{prop:li}.
\end{proof}

Suppose now that $\cM=\cF(\mathbf{E},\mu)$ for  an admissible $n$-coalgebra bundle  $(\mathbf{E}\to M,\mu)$.
Recall that, by Theorem \ref{thm:compatder}, 
$$
\vf^{-k}=\Gamma_{\der({\bf E}^*_\bR, m)_{-k}}\quad \text{for}\quad  k=0,\ldots, n.
$$  

By the previous lemma, there are vector subbundles $(C_k\to M) \subseteq (\der({\bf E}^*_\bR, m)\to M)$ such that 
\begin{equation}\label{eq:ck}
    \Gamma_{C_k} = \cD\cap \vf^{-k}, \qquad  k=0, \ldots, n.
\end{equation}

Recall the maps $\Theta_{i,j}$ in \eqref{eq:theta} and $\mathrm{res}_r$ (see Prop.~\ref{prop:pullbacks}).

\begin{theorem}\label{thm:dist}
    Let $(\mathbf{E}\to M,\mu)$ be an admissible $n$-coalgebra bundle. A distribution in $\cM=\cF(\mathbf{E},\mu)$ is equivalent to
    a collection of vector subbundles $(C_k\to M)\subseteq \der(\mathbf{E}^*_\mathbb{R}, m)_{-k}$, for $k= 0,\cdots,n$, satisfying
    \begin{equation}\label{eq:distcond}
    \Theta_{1,-(k+1)}(E_{-1}^*\otimes C_{k+1}) + \ldots + \Theta_{n-k,-n}(E_{k-n}^* \otimes C_{n})=\ker(\mathrm{res}_k) \cap C_k.
    \end{equation}
Moreover, for such a collection $\{C_k\}$, the intersections $\ker(\mathrm{res}_k) \cap C_k$ have constant rank, and local generators of the distribution are identified with the union of local frames of $C_k/(\ker(\mathrm{res}_k) \cap C_k)$, for $k=0,\ldots,n$.  (Note that \eqref{eq:distcond} is trivial for $k=n$ since $\mathrm{res}_n$ is the identity map.)
\end{theorem}
\begin{proof}
Given a distribution $\cD$ in $\cM$, condition \eqref{eq:distcond} is a direct consequence of \eqref{eq:ck} and \eqref{eq:localgen}.

On the other hand,  consider subbundles $C_k \subseteq \der(\mathbf{E}^*_\mathbb{R}, m)_{-k}$, for $k= 0,\cdots,n$, satisfying \eqref{eq:distcond}. Using that $\Gamma_{C_k} \subseteq \Gamma_{\der({\bf E}^*_\bR, m)_{-k}} = \vf^{-k}$, we consider the subsheaf of $C_\cM$-modules of $\vf$ given by
$$
\cD = C_\cM \cdot \Gamma_{C_0} + \ldots + C_\cM \cdot \Gamma_{C_n}.
$$
We must check that it is a distribution, in the sense of Definition~\ref{def:dist}.

Note that, as a sheaf of $C^\infty_M$-modules,  $C_\cM^1\cdot \Gamma_{C_{k+1}}+\ldots + C_\cM^{n-k}\cdot \Gamma_{C_n}$ is locally free. By \eqref{eq:distcond}, $\ker(\mathrm{res}_k) \cap C_k$ is a vector subbundle of $C_k$, so that
$$
C_\cM^1\cdot \Gamma_{C_{k+1}}+\ldots + C_\cM^{n-k}\cdot \Gamma_{C_n} = \Gamma_{\ker(\mathrm{res}_k) \cap C_k}.
$$
In particular, the fact that $C_\cM^1\cdot \Gamma_{C_{k+1}}+\ldots + C_\cM^{n-k}\cdot \Gamma_{C_n}\subseteq \Gamma_{C_k}$ implies that
$$
 \cD\cap \vf^{-k} = \Gamma_{C_k}. 
$$
With the choice of splittings, we may write $C_k = (\ker(\mathrm{res}_k) \cap C_k) \oplus ({C_k}/({\ker(\mathrm{res}_k) \cap C_k}))$,
and consider local frames of ${C_k}/({\ker(\mathrm{res}_k) \cap C_k})$ as given by local vector fields $\{X_k^{i_k}\}$, $i_k=1,\ldots,d_k$, of degree $-k$. The union of these vector fields for $k=0,\ldots,n$ on any small open subset $U\subseteq M$ generates $\cD(U)$, as in Definition~\ref{def:dist}.

\end{proof}

\begin{remark}[Restrictions] \label{rem:rest}
For an $n$-manifold $\cM$, recall the tower of $\mathbb{N}$-manifolds,
$$
M=\cM_0\leftarrow\cM_1\cdots\leftarrow\cM_{n-1}\leftarrow\cM_{n}=\cM.
$$
Given a distribution $\cD$ on $\cM$ of rank $d_0|\cdots|d_n$, there is a distribution $\cD_r$ on each $\cM_r$, of rank $d_0|\cdots|d_r$, characterized by the condition
$$
\cD_r\cap \cT^{-k}_{\cM_r} = \mathrm{res}_r(\cD \cap \vf^{-k}),
$$
for $k=0,\ldots,r$. This can be shown directly using local generators, or the viewpoint of Thm.~\ref{thm:dist}. Moreover, if $\cD$ is involutive so is each $\cD_r$.  
\hfill $\diamond$
\end{remark}

\begin{example}[Distributions on $1$-manifolds]\label{1-dist} 
Any $1$-manifold is of the form $E[1]$ for a vector bundle $E\to M$, and, as seen in Example \ref{1-vf},  
$$
\cT_{E[1]}^0=\Gamma_{\mathrm{Der}(E^*)} \quad \mbox{ and } \quad  \cT_{E[1]}^{-1}=\Gamma_E.
$$
It follows from Theorem~\ref{thm:dist} that distributions on $E[1]$ are described by vector subbundles
 $C_0\subseteq \der(E^*)$ and $C_{1}\subseteq E$  satisfying 
$$
C_{1}\otimes E^* = C_0\cap \ker(\mathrm{res}_0) =   C_0\cap \ker\sigma,
$$
where $\sigma: \der(E^*)\to TM$ is the symbol map
(see Remark~\ref{rem:symbol}). Hence, for $n=1$, we recover the geometric characterization of distributions in \cite[Lemma~3.1]{zam:dis}. 

A distribution represented by $C_0$ and $C_1$ is involutive if and only  if $C_0$ is a Lie subalgebroid of $\der(E^*)\cong \der(E)$ and $C_{1}\subseteq E$ is invariant under $C_0$. In this case, as shown in 
\cite[Lemma 3.2]{zam:dis}, $C_0$ is equivalent to the data of an involutive distribution $F\subseteq TM$ and a flat $F$-connection $\nabla$ on $E/C_{1}$: indeed, given $C_0$, we set  
$F=\sigma(C_0)$, and
$$
\nabla_X[e]=[D^{\dagger}( e)] \quad \text{for}\ X\in\Gamma_F, \ [e]\in\Gamma_{E/C_{1}},
$$ 
and $D^\dagger$ the dual of $D\in\Gamma({C_0})$ with $\sigma(D)=X$; given $F$ and $\nabla$, we define $C_0$ such that
\begin{equation}\label{eq:C0}
\Gamma_{C_0}=\{ D \in \Gamma_{\der(E^*)}\,|\, \sigma(D) \in \Gamma_F, \, D^\dagger(\Gamma_{C_{1}})\subseteq \Gamma_{C_{1}}, \, [D^\dagger]=\nabla_{\sigma(D)} \},
\end{equation}
where $[D^\dagger] \in \Gamma_{\der(E/C_{1})}$ is given by $[D^\dagger][e]=[D^\dagger (e)]$, for $[e]$ denoting the class of $e\in \Gamma_E$ in $\Gamma_{E/C_{1}}$.
\hfill $\diamond$
\end{example}

The next example is a particular instance of the previous one.

\begin{example}[Lie algebroids] 
A Lie algebroid structure on $A\to M$, with anchor $\rho$ and bracket $[\cdot,\cdot]$, is equivalent to  a degree $-1$ Poisson structure $\Pi$ on the manifold $A^*[1]$ that can be written in local coordinates $\{ x^i, e_\alpha\}$  as
$$
\Pi=\rho^i_\alpha(x) \frac{\partial}{\partial x^i}\wedge\frac{\partial}{\partial e_\alpha}+ c^{\gamma}_{\beta\alpha}(x) e_{\gamma}\frac{\partial}{\partial e_\beta}\wedge\frac{\partial}{\partial e_\alpha}.
$$
We have the subsheaf of hamiltonian vector fields $\cA_{Ham}\subseteq \cT_{A^*[1]}$, defined on open subsets $U\subseteq M$ by 
$$
\cA_{Ham}(U)=\{ \Pi(df)\in\cT_{A^*[1]}(U)\ |\ f\in C_{A^*[1]}(U)\},
$$
that is locally generated by the vector fields
$$
 X^\alpha=\rho^i_\alpha(x) \frac{\partial}{\partial x^i}+ c^{\gamma}_{\beta\alpha}(x) e_{\gamma}\frac{\partial}{\partial e_\beta}, \quad \mbox{and} \;\;\; Y^i=\rho^i_\alpha(x)\frac{\partial}{\partial e_\alpha}.
 $$
Then $\cA_{Ham}$ defines a distribution (in the sense of Definition~\ref{def:dist}) if and only if $\rho:A\to TM$ has constant rank and $\ker(\rho)$ is a bundle of abelian Lie algebras. The corresponding triple $(C_{1},F,\nabla)$ as in Example \ref{1-dist}
is given by 
$C_{1}=\im(\rho^t)\subseteq A^*$, $F=\im\rho$,  and $\nabla$ the dual of the $F$-connection given by 
$$
\nabla^\dagger:\Gamma_{\im(\rho)}\times\Gamma_{\ker(\rho)}\to\Gamma_{\ker(\rho)}, \qquad \nabla^\dagger_Xe=[\widehat{e}, e], 
$$
where $\widehat{e} \in \Gamma_A$ is such that $\rho(\widehat{e})=X$, noticing that $\ker(\rho)^*= A^*/\im(\rho^t)= A^*/C_1$.
\hfill $\diamond$

\end{example}

A source of examples of $\mathbb{N}$-manifolds carrying distributions in degree $2$ is given by coisotropic submanifolds of symplectic 2-manifolds, as discussed in \cite[Sec.~4]{bur:super}.

It follows from Example~\ref{1-dist} that a distribution on a 1-manifold $E[1]$ of rank $d_0|0$ is specified by an involutive subbundle $F\subseteq TM$ and a flat $F$-connection on $E$ (or, dually, on $E^*$). For later use, we will need the extension of this result to $\mathbb{N}$-manifolds of arbitrary degrees.

\begin{proposition}\label{prop:0-dist}
     Let $(\mathbf{E}\to M, \mu)$ be an admissible $n$-coalgebra bundle and $\cM=\cF(\mathbf{E}, \mu)$ the associated $n$-manifold. An involutive distribution $\cD$ on $\cM$ of rank $d_0|0|\cdots|0$ is equivalent to the following data: an involutive subundle $F\subseteq TM$ and flat $F$-connections $\nabla^i$ on $E_{-i}^*$, $i=1,\dots,n$, satisfying 
     \begin{equation}\label{eq:nablacomp}
     \nabla^{i+j}_X(\mu^*(e, e'))=\mu^*(\nabla^i_Xe, e')+\mu^*(e, \nabla_X^j e')
     \end{equation}
for all $e\in\Gamma({E_{-i}^*})$, $e'\in\Gamma({E_{-j}^*})$, and $X\in \Gamma(F)$.     
\end{proposition}
\begin{proof}
According to Theorem~\ref{thm:dist}, a distribution $\cD$ with rank $d_0|0|\ldots|0$ is equivalent to a vector subbundle $C_0 \subseteq \der(\mathbf{E}_\bR^*,m)_0$ such that $\ker(\sigma) \cap C_0=0$. Recall that sections of $\der(\mathbf{E}_\bR^*,m)_0$ are given by derivations $D^i$ of $E_{-i}^*$, all with the same symbol $X=\sigma(D^i)$, 
such that
$$
D^{i+j}(\mu^*(e, e'))=\mu^*(D^ie, e')+\mu^*(e, D^j e').
$$
The symbol map $\sigma$ projects $C_0$ isomorphically onto its image $F:=\sigma(C_0)$, so $C_0$ is the image of $\nabla = (\sigma|_{C_0})^{-1}: F \to \der(\mathbf{E}_\bR^*,m)_0$. The distribution $\cD$ is involutive if and only if $C_0$ is a subalgebroid, which amounts to $F$ being involutive and $\nabla$
being determined by $F$-connections $\nabla^i$ on $E_{-i}^*$, satisfying \eqref{eq:nablacomp}, which are flat. 
\end{proof}

\section{The Frobenius theorem}\label{sec:frobenius}

In this section we will prove the Frobenius theorem for $n$-manifolds, asserting that any involutive distribution is locally generated by coordinate vector fields.

As a warm up, we start by proving the Frobenius theorem for $1$-manifolds.
Let $E\to M$ be a vector bundle with $m_0=\dim(M)$ and $m_1=\mathrm{rank}(E)$.

\begin{proposition} \label{Frob1}
 Let $\cD$ be an involutive distribution on $E[1]$ of rank $d_0|d_1$. Then any point in $M$ admits a neighborhood $U$ with coordinates $\{x^\alpha, e^{ \beta}\}$, $\alpha=1,\ldots, m_0$ and $\beta=1,\ldots,m_1$, on $E[1]_{|U}$ such that
\begin{equation}\label{eq:frob1}
 \cD(U)=\langle \frac{\partial}{\partial x^1},\cdots,\frac{\partial}{\partial x^{d_0}},\frac{\partial}{\partial e^1},\cdots, \frac{\partial}{\partial e^{d_1}}\rangle.
\end{equation}
\end{proposition}

\begin{proof}
Any point in $M$ admits a neighborhood $U$ such that  
\begin{equation}\label{eq:gener}
\cD(U)=\langle X^1, \cdots, X^{d_0}, Y^1, \cdots, Y^{d_1}\rangle,
\end{equation}
where $X^i\in\cT^0_{E[1]}(U), \ Y^j\in\cT^{-1}_{E[1]}(U)$ are linearly independent. As recalled in  Example \ref{1-dist}, $\cD$ is described by a triple $(C_1, F, \nabla)$, where $C_{1}\subseteq E$ is a subbundle, $F\subseteq TM$ is an involutive subbundle and $\nabla$ is a flat $F$-connection on $E/C_1$. With respect to these data, picking local generators of $\cD$ as in \eqref{eq:gener} amounts to choosing a frame
$\{\xi_1, \cdots \xi_{d_1}\}$ of ${C_1}_{|_U}$  and local sections 
$\{D_1,\ldots,D_{d_0}\}$ of ${C_0}_{|_U}$, following the notation in \eqref{eq:C0}, so that the symbols $\sigma(D_j)$, $j=1,\ldots,d_0$, form a frame for
$F_{|_U}$. 

Showing that one can find local generators of $\cD$ as in \eqref{eq:frob1} translates into showing that there exist local coordinates $\{x^\alpha\}_{\alpha=1}^{m_0}$ on $U$ and a frame $\{e^\beta\}_{\beta=1}^{m_1}$ of $E^*_{|U}$, as well as a frame 
$\{\xi_1, \cdots \xi_{d_1}\}$ of ${C_1}_{|_U}$ and sections
$\{D_1,\ldots,D_{d_0}\}$ of ${C_0}_{|_U}$ as above  satisfying
\begin{equation}\label{eq:frobconds}
\sigma(D_a)(x^\alpha)=\delta_{a\alpha}, \;  \quad D_a(e^\beta)=0, \quad \xi_b (e^\beta)=\delta_{b \beta},
\end{equation}
for $a= 1,\ldots, d_0$, $\alpha=1,\ldots, m_0$, $b= 1,\ldots, d_1$, $\beta = 1, \ldots, m_1$.

We construct the above data $\{x^\alpha\}$, $\{e^\beta\}$, $\{\xi_b\}$ and $\{D_a\}$ as follows. Shrinking $U$ if necessary, the usual Frobenius theorem implies the existence of local coordinates $\{x^\alpha\}_{\alpha=1}^{m_0}$  on $U$ such that\footnote{For local coordinates $\{x^\alpha\}$ on $U\subseteq M$, we denote the corresponding coordinate vector fields on $U$ by $\partial_{x^\alpha}$; when $\{x^\alpha\}$ is part of a  coordinate system $\{x^\alpha, e^\beta\}$ on $\cM_{|_U} =  {E[1]}_{|_U}$, we use the notation $\frac{\partial}{\partial x^\alpha}$ for the corresponding degree $0$ vector field (which is a derivation of ${E^*}_{|_U}$ with symbol $\partial_{x^\alpha}$).} $F_{|_U}=\langle {\partial_{x^a}} \ | \ a=1,\ldots, d_0\rangle$,  we can assume that $E_{|U}$ and $C_{1|U}$ are trivial and,  writing $E_{|U}= C_{1}|_{U} \oplus (E/C_{1})_{|U}$ (with the choice of a splitting), we define a flat $F$-connection on $E_{|U}$ by 
$$
\widehat{\nabla}=\nabla^{triv}\oplus\nabla,
$$
where $\nabla^{triv}$ denotes the trivial $F$-connection on $C_{1|U}$.
Let $\widehat{\nabla}^{\dagger}$ be the dual flat $F$-connection on $E^*_{|U}$, and let $\{e^\beta\}_{\beta=1}^{m_1}$ be a frame of flat sections of $E^*_{|U}$ such that $\{ e^1,\cdots e^{d_1}\}$  is a frame of flat section of $C^*_{1}|_{U}$.  Take $\{\xi_b\}_{b=1}^{d_1}$ to be the dual frame on $C_{1}|_{U}$ and
$$
D_a:= \widehat{\nabla}^{\dagger}_{{\partial_{x^a}}},\qquad a=1,\cdots, d_0.
$$
These data satisfy \eqref{eq:frobconds}.
 \end{proof}

We will now generalize the previous proposition to $n$-manifolds of arbitrary degree; for that, we will first consider the special cases of distributions of ranks $d_0|0|\ldots|0$ and $0|d_1|\ldots|d_n$.

\begin{lemma}\label{L3}
Let $\cM$ be an $n$-manifold of dimension $m_0|\ldots|m_n$ with an involutive distribution   $\cD$ of rank $d_0|0|\cdots|0$. Then any point in $M$ admits a neighborhood $U$ with coordinates $\{x^\alpha, e^{ \beta_j}_j\}$ on $\cM_{|U}$ such that
\begin{equation}\label{eq:xx}
 \cD(U)=\langle \frac{\partial}{\partial x^1},\cdots, \frac{\partial}{\partial x^{d_0}}\rangle.
\end{equation}  
\end{lemma}

\begin{proof}
    Let $(\mathbf{E}\to M,\mu)$ be an $n$-coalgebra bundle such that $\cM=\cF(\mathbf{E},\mu)$. As seen in Proposition \ref{prop:0-dist}, the distribution $\cD$ is equivalent to an involutive subbundle $F\subseteq TM$ and flat  $F$-connections $\nabla^j$ on $E^*_{-j}$, for $j=1,\ldots, n$, satisfying \eqref{eq:nablacomp}; with respect to this geometric description, the statement of the lemma translates into the existence of a neighborhood $U$ of any point in $M$ admitting coordinates $\{x^\alpha\}_{\alpha=1}^{m_0}$, together with a frame $\{Z_1,\ldots,Z_{d_0}\}$ of $F_{|_U}$ and  frames $\{ e_j^{\beta_j}\}_{\beta_j=1}^{m_j}$ of $(\ker \mu_{-j})^*$, for $j=1,\ldots,n$,  satisfying the following properties: 
    $$
    Z_a(x^\alpha)=\delta_{a\alpha}, \qquad \nabla^j_{Z_a}(e_j^{\beta_j}) =0. 
    $$
    In the equation on the right-hand side, we regard $e_j^{\beta_j}$, $\beta_j=1,\ldots,m_j$, as sections of ${E_{-j}^*}_{|_U}$ by means of suitable splittings of the exact sequences
$$
0\to (K_{-j}^\mu)^*\xrightarrow{} E_{-j}^*\to \ker(\mu_{-j})^*\to 0,
$$
for $j=1,\ldots,n$, see Remark~\ref{rem:local}.

The desired data $\{x^\alpha\}$, $\{Z_a\}$, $\{e^{\beta_j}_j\}$ is constructed as follows. Taking $U$  small enough, the usual Frobenius theorem implies the existence of local coordinates $\{x^\alpha\}_{\alpha=1}^{m_0}$  on $U$ such that $\{Z_a:=\partial_{x_a}\}_{a=1}^{d_0}$ is a frame for 
$F_{|_U}$. 
It follows from condition \eqref{eq:nablacomp} that $(K_{-j}^\mu)^* \subseteq E_{-j}^*$ is $\nabla^j$-invariant, and this ensures the existence of linearly independent flat sections
$\{e^{\beta_j}_j\}_{\beta_j=1}^{m_j}$ of ${E_{-j}^*}_{|_U}$ spanning a subbundle complementary to $(K_{-j}^\mu)^*$.
\end{proof}

\begin{lemma}\label{L2}
  Let $\cM$ be an $n$-manifold of dimension $m_0|\ldots|m_n$ and  $\cD$ an involutive distribution of rank $0|d_1|\cdots|d_n$. Then any point in $M$ admits a neighborhood $U$ with coordinates $\{x^\alpha, e^{ \beta_j}_j\}$ on $\cM_{|U}$ such that
\begin{equation*}
 \cD(U)=\langle \frac{\partial}{\partial e_1^1},\cdots, \frac{\partial}{\partial e_1^{d_1}}, \ldots,  \frac{\partial}{\partial e_n^1},\cdots, \frac{\partial}{\partial e_n^{d_n}}\rangle.
\end{equation*}  
\end{lemma}
\begin{proof}
We proceed by induction on $n$. For $n=1$ the statement is a particular case of Prop.~ \ref{Frob1}. Assuming that the lemma holds for $\mathbb{N}$-manifolds of degree $n-1$, we will prove it for degree $n$. For simplicity, we will use the notation 
$$
p^\gamma=e^{\gamma}_n
$$ 
for degree $n$ coordinates. 

Following Remark~\ref{rem:rest}, $\cD$ defines a distribution $\cD_{n-1}$ on the truncation $\cM_{n-1}$. The explicit local picture is as follows: in an open $U$ with coordinates $\{x^\alpha, e_j^{\beta_j}, p^\gamma\}$, for $j=1,\ldots, n-1$, on $\cM_{|_U}$, and such that $\cD(U)$ is generated by linearly independent vector fields ${Y}_j^{b_j}$ and $Z^c$, where $j=1,\ldots,n-1$, $b_j=1,\ldots, d_j$, $c=1,\ldots, d_n$, $Y^{b_j}_j$ has degree $-j$ and $Z^c$ has degree $-n$, the truncation ${\cM_{n-1}}_{|_U}$ has coordinates $\{x^\alpha, e_j^{\beta_j}\}$ and $\cD_{n-1}(U)$ is generated by the restrictions of ${Y}_j^{b_j}$ to ${\cM_{n-1}}_{|_U}$.

By the induction hypothesis, shrinking $U$ if necessary, we may assume that
\begin{equation*}
 \cD_{n-1}(U)=\langle \frac{\partial}{\partial e_j^1},\cdots, \frac{\partial}{\partial e_j^{d_j}}\ |\ j=1, \ldots, n-1\rangle,
\end{equation*} 
and that the vector fields $Y_j^{b_j}$ and $Z^c$ satisfy
\begin{equation}\label{trunX}
  \mathrm{res}_{n-1}({Y_j}^{b_j})=\frac{\partial}{\partial e_j^{b_j}}\quad \text{and}\quad \mathrm{res}_{n-1}(Z^c)=0.    
\end{equation}

The vector fields $Z^c$ form a local frame for the subbundle  $C_{n}\subseteq\ker \mu_{-n}\subseteq E_{-n}$ (see Theorems~\ref{thm:compatder} and \ref{thm:dist}). Passing to a smaller neighborhood if necessary, we can extend this frame to a local frame of $\ker\mu_{-n}$, and the dual frame of $(\ker\mu_{-n})^*$
defines degree $n$ coordinates 
$$
\{p^1, \ldots p^{d_n}, \widehat{p}^{d_n+1}, \ldots, \widehat{p}^{m_n} \}
$$ 
that complete $\{x^\alpha,e_j^{\beta_j}\}$ to a set of local coordinates of $\cM_{|U}$ in such a way that
$Z^c=\frac{\partial}{\partial p^c}$
and, by \eqref{trunX},
$$
{Y}^{b_j}_j=\frac{\partial}{\partial e^{b_j}_j}+\sum_{l=1}^{d_n}{Y}^{b_j}_j(p^l)\frac{\partial}{\partial p^l}+\sum_{l=d_n+1}^{m_n}{Y}^{b_j}_j(\widehat{p}^l)\frac{\partial}{\partial \widehat{p}^l},
$$
cf. Prop.~\ref{neg-gen-vf}. By considering linear combinations, we
may redefine 
$$
{Y}_j^{b_j} :=\frac{\partial}{\partial e^{b_j}_j}+\sum_{l=d_n+1}^{m_n} {Y}_j^{b_j}(\widehat{p}^l) \frac{\partial}{\partial \widehat{p}^l}
$$
and conclude that 
$$
\cD(U)=\langle  {Y}_j^{b_j},\; \frac{\partial}{\partial p^c}  \rangle,
$$
for $j= 1,\ldots n-1$, $b_j= 1, \ldots, d_j,$ and $c= 1, \ldots, d_n$. 
We will now implement successive changes of coordinates just involving the variables $\widehat{p}^l$ that will turn each $Y_j^{b_j}$ into $\frac{\partial}{\partial e_j^{b_j}}$.

We will simplify the notation and denote the coordinates $e^{b_j}_j$ by $e^b$.
Suppose that we have already brought $s-1$ vector fields (note that $s-1$ could be zero) to the desired form, so we are in the situation
\begin{equation}\label{eq:DU}
\cD(U)=\langle \frac{\partial}{\partial e^{1}}, \ldots, \frac{\partial}{\partial e^{s-1}},\ {Y}^s, \ldots,  {Y}^k, \frac{\partial}{\partial p^1},\ldots, \frac{\partial}{\partial p^{d_n}} \rangle,
\end{equation}
where ${Y}^b=\frac{\partial}{\partial e^b}+\sum_{l=d_n+1}^{m_n} {g}^{bl}\frac{\partial}{\partial \widehat{p}^l}$. Note that, since each coefficient $g^{bl}$ has degree $n-|e^b|<n$, we have that
\begin{equation}\label{eq:gij}
\frac{\partial}{\partial p^c} g^{bl}=0, \qquad \frac{\partial}{\partial \widehat{p}^i} g^{bl}=0,
\end{equation}
for $c=1,\ldots,d_n$ and $i=d_n+1,\ldots,m_n$.

We will now define a change of variables turning $Y^s$ into $\frac{\partial}{\partial e^s}$.
The following is a key observation.

\smallskip

\noindent {\bf Claim}. {\em For each $l$ with $g^{sl}\neq 0$, there is a (unique) degree $n$ function $G^{sl}$ on $\cM_{|U}$ such that $\frac{\partial}{\partial e^s}G^{sl}= g^{sl}$.}

\smallskip

Let us verify this claim.
Since $g^{sl}$ depends  polynomially on $e^s$, we may write it as 
\begin{equation}\label{eq:as}
g^{sl}= \sum_{t=0}^r (e^s)^t a^{sl}_t
\end{equation}
(note that this sum must have finitely many terms), and we simply set 
\begin{equation}\label{eq:Gs}
G^{sl} = \sum_{t= 0}^r \frac{1}{1+t} (e^s)^{t+1} a^{sl}_t.
\end{equation} 
Such $G^{sj}$ satisfies the condition in the claim provided it is nonzero. This is clearly the case when $|e^s|$ is even. To check that this also holds when $|e^s|$ is odd, we will make use of the involutivity condition of $\cD$. Using the second condition in \eqref{eq:gij}, we see that (recall that $|Y^s|=-|e^s|$ is odd)
$$
[Y^s,Y^s] = 2 \sum_{l=d_n+1}^{m_n}  \frac{\partial g^{sl}}{\partial e^s} \frac{\partial}{\partial \widehat{p}^l},
$$
and this vector can only be in $\cD(U)$ if 
\begin{equation}\label{eq:invol1}
\frac{\partial g^{sl}}{\partial e^s} = 0.
\end{equation}
Note that $G^{sl}= e^s g^{sl}$ when $|e^s|$ is odd, and \eqref{eq:invol1} ensures that this is nonzero. This concludes the proof of the claim.

Consider the change of coordinates that only modifies the coordinates $\widehat{p}^l$ according to 
$$
\widehat{p}^l \mapsto \widehat{p}^l - G^{sl}, \qquad l=d_n+1 ,\ldots, m_n,
$$
while keeping all the other coordinates unchanged (here we set $G^{sl}$ to be zero if $g^{sl}$ is zero). 

To verify how this transforms the generators of $\cD(U)$ in \eqref{eq:DU}, we first observe that
\begin{equation}\label{eq:invol2}
\frac{\partial g^{sl}}{\partial e^b} =0, \qquad b=1,\ldots, s-1,
\end{equation}
which is a consequence of
$$
[\frac{\partial}{\partial e^b}, Y^s] = \sum_{l=d_n+1}^{m_n} \frac{\partial g^{sl}}{\partial e^b} \frac{\partial}{\partial \widehat{p}^l},
$$
and the involutivity of $\cD$. Note that conditions \eqref{eq:gij} and \eqref{eq:invol2}
hold for $a_t^{sl}$ in \eqref{eq:as}, and hence (see \eqref{eq:Gs})
\begin{equation}\label{eq:Gsj}
\frac{\partial G^{sl}}{\partial p^c} =0, \quad \frac{\partial G^{sl}}{\partial \widehat{p}^i} =0, \; \mbox{ and } \; \frac{\partial G^{sl}}{\partial e^b} =0,
\end{equation}
for $b=1,\ldots, s-1$. It then follows, by the usual transformation rule for vector fields under change of coordinates, that
$$
\frac{\partial}{\partial p^c}, \quad \frac{\partial}{\partial \widehat{p}^i}, \;\; \mbox{ and } \;\; \frac{\partial}{\partial e^b},
$$
for $b=1,\ldots, s-1$, remain unchanged. 
For $b\geq s$ we have that
$\frac{\partial}{\partial e^b}$ transforms to 
$$
\frac{\partial}{\partial e^b} -  \sum_{l=d_n+1}^{m_n} \frac{\partial G^{sl}}{\partial e^b} \frac{\partial}{\partial \widehat{p}^l}
$$ 
in the new coordinates. In particular, the fact that $\frac{\partial}{\partial e^s}G^{sl}= g^{sl}$ implies that
$Y^s$ agrees with $\frac{\partial}{\partial e^s}$ after the change of coordinates.

We conclude that, in the new coordinates, we can write
$$
\cD(U)=\langle \frac{\partial}{\partial e^{1}}, \ldots, \frac{\partial}{\partial e^{s}},\ {Y}^{s+1}, \ldots,  {Y}^k, \frac{\partial}{\partial p^1}, \ldots, \frac{\partial}{\partial p^{d_n}}  \rangle,
$$
where ${Y}^b=\frac{\partial}{\partial e^b}+\sum_{l=d_n+1}^{m_n} {\tilde{g}}^{bl}\frac{\partial}{\partial \widehat{p}^l}$. The proof of the lemma follows from a successive repetition of the argument. 
\end{proof}

We now prove our main theorem.

\begin{theorem}[Frobenius theorem for $\mathbb{N}$-manifolds]\label{localFrob}
Let $\cM$ be an $n$-manifold of dimension $m_0|\ldots|m_n$ and $\cD$ an involutive distribution on $\cM$ of rank $d_0|\cdots|d_n$. Then any point in $M$ has a neighborhood $U$ with coordinates $\{x^\alpha, e_j^{ \beta_j}\}$ on $\cM_{|U}$ such that
\begin{equation*}
 \cD(U)=\langle \frac{\partial}{\partial x^1},\ldots,\frac{\partial}{\partial x^{d_0}},\frac{\partial}{\partial e_1^{1}},\ldots, \frac{\partial}{\partial e_1^{d_1}},\ldots, \frac{\partial}{\partial e_n^{1}},\ldots, \frac{\partial}{\partial e_n^{d_n}}\rangle.
\end{equation*}
\end{theorem}

\begin{proof}
Let $U$ be a neighborhood where $\cD(U)$ is generated by vector fields $\{X_k^{i_k}\}$, $k=0,\ldots,n$, as in Definition~\ref{def:dist}. By applying Lemma~\ref{L2} to the distribution generated by $\{X_k^{i_k}\}$ for $k>0$, we conclude that (by shrinking $U$ if necessary) that there exist coordinates $\{{x}^\alpha, e_j^{ \beta_j}\}$ on $\cM_{|U}$ such that
\begin{equation*}
 \cD(U)=\langle \widehat{X}^1,\cdots, \widehat{X}^{d_0},\frac{\partial}{\partial e_1^{1}},\cdots, \frac{\partial}{\partial e_1^{d_1}},\cdots, \frac{\partial}{\partial e_n^{1}},\cdots, \frac{\partial}{\partial e_n^{d_n}}\rangle,
\end{equation*}
where $\widehat{X}^1,\cdots,\widehat{X}^{d_0}$ are linearly independent, degree $0$ vector fields. In the above coordinates,
$$
\widehat{X}^a = \sum_{\alpha=0}^{m_0} g^a_\alpha\frac{\partial}{\partial x^\alpha} + \sum_{j=1}^n\sum_{\beta_j=1}^{m_j} f^{aj}_{\beta_j} \frac{\partial}{\partial e_j^{\beta_j}}
$$
with $g^a_\alpha\in C^\infty(U)=C^0_{\cM}(U)$ and $f^{aj}_{\beta_j}\in C^j_\cM(U).$  Hence, by taking linear combinations, we can write
\begin{equation*}
 \cD(U)=\langle {X}^1,\cdots, {X}^{d_0},\frac{\partial}{\partial e_1^{1}},\cdots, \frac{\partial}{\partial e_1^{d_1}},\cdots, \frac{\partial}{\partial e_n^{1}},\cdots, \frac{\partial}{\partial e_n^{d_n}}\rangle,
\end{equation*}
where $X^a$ is of the form
$$
 X^a=\sum_{\alpha=0}^{m_0} g^a_\alpha\frac{\partial}{\partial x^a}+\sum_{j=1}^n\sum_{b_j=d_j+1}^{m_j} f^{aj}_{b_j}\frac{\partial}{\partial e_j^{b_j}}.
$$
Note that, for $i=1,\ldots,n$ and $\beta_i=1,\ldots,d_i$,
$$
[\frac{\partial}{\partial e_i^{\beta_i}}, X^a]= \sum_{j=1}^n\sum_{b_j=d_j+1}^{m_j}\frac{\partial f^{aj}_{b_j}}{\partial e_i^{\beta_i}}\frac{\partial}{\partial e_j^{b_j}}
$$
so it follows from the involutivity of $\cD$ that
\begin{equation}\label{eq:fdep}
\frac{\partial f^{aj}_{b_j}}{\partial e_i^{\beta_i}}=0.
\end{equation}

By splitting the local coordinates on $\cM_{|U}$ as $\{x^\alpha, e_j^{b_j}\}$ and $\{e_j^{\beta_j}\}$, for $j=1,\ldots,n$, $b_j=d_j+1,\ldots,m_j$, and $\beta_j=1,\ldots,d_j$  we write
$$
\cM_{|U} = \bR^{m_0|m_1-d_1|\cdots|m_n-d_n}_{|U} \times {\bR^{0|d_1|\cdots|d_n}}
$$ 
in such a way that, by \eqref{eq:fdep},  $\cD$ is identified with the direct product of the involutive distributions 
$$
\cD_1 =\langle X^1,\ldots X^{d_0}\rangle
$$ 
on $\bR^{m_0|m_1-d_1|\cdots|m_n-d_n}_{|U}$ of rank $d_0|0|\cdots|0$ and $\cD_2 =\langle \frac{\partial}{\partial e_j^1},\cdots, \frac{\partial}{\partial e_j^{d_j}}\ |\ j=1,\ldots,n \rangle$
on $\bR^{0|d_1|\cdots|d_n}$. We conclude the proof by applying Lemma \ref{L3} to $\cD_1$.

\end{proof}

The following is a simple instance of the Frobenius theorem. 

\begin{corollary}\label{cor:vf}
Let $X$ be a vector field of degree $-k$ on an $n$-manifold $\cM$, for $k=1,\ldots, n$ (resp. $k=0$), such that  $X_x\neq 0$ for $x \in M$ and $[X,X]=\lambda X$, where $\lambda$ is a degree $k$ function on $\cM$. Then there is a neighborhood of $x$ with coordinates $\{x^\alpha, e_j^{\beta_j}\}$ such that $X=\frac{\partial}{\partial e_k^1}$ (resp. $X=\frac{\partial}{\partial x^1}$).
\end{corollary}

This last result follows from the observation that the vector field $X$ generates an involutive distribution in a neighborhood of $x$. Analogous local normal forms for vector fields on supermanifolds can be found in \cite{sha:vec} (cf. \cite[Lemmas~4.7.5, 4.7.6]{var:book}), see also \cite{mon:exis}.

\begin{remark}[Normal forms of $Q$-structures] 
Since vector fields of positive degree vanish when evaluated at points in the body, they do not generate distributions and hence do not fit into the framework of the Frobenius theorem. For homological vector fields (see Example~\ref{ex:Qman}), the study of local normal forms is motivated by several results in geometry; e.g., local splitting theorems for Lie algebroids as well as similar results for Courant and higher Lie algebroids  (see \cite{bur:def}) can be viewed as local normal forms for $Q$-structures on $\mathbb{N}$-manifolds (see \cite{vai:vec} for results on supermanifolds). 
\hfill $\diamond$
\end{remark}

\begin{remark}[Integrability and foliations]
We say that a  distribution $\cD$ of rank $d_0|\ldots|d_n$ on an $n$-manifold $\cM$ of dimension $m_0|\ldots|m_n$ is {\em integrable} if there exist, around any point in $M$, a neighborhood $U$, neighborhoods $V_1$ and $V_2$ of the origin in  $\bR^{d_0}$ and ${\bR^{m_0-d_0}}$, respectively, and an isomorphism 
\begin{equation}\label{eq:folchart}
\Psi:\cM_{|U}\to{\bR^{d_0|\cdots|d_n}}_{|{V_1}}\times {\bR^{m_0-d_0|\cdots|m_n-d_n}}_{|V_2}
\end{equation}
identifying $\cD$ with the distribution $\cT_{\bR^{d_0|\cdots|d_n}} \times 0$. As in the classical case, Theorem~\ref{localFrob} is equivalent to saying that a distribution $\cD$ is integrable if and only if it is involutive. 

A {\em foliation atlas} on $\cM$ is a collection of charts as in \eqref{eq:folchart} such that the transition maps preserve the distribution $\cT_{\bR^{d_0|\cdots|d_n}} \times 0$. As in classical geometry (see e.g. \cite{moe:lie}), one can define foliations on $\bN$-manifolds by means of foliation atlases and establish an equivalence between foliations and involutive distributions.
\hfill $\diamond$
\end{remark}

\end{document}